\documentclass[11pt]{article}
\usepackage{amsmath}
\usepackage{amssymb}
\usepackage{latexsym}

\newtheorem{df}{Definition}[section]
\newtheorem{thm}[df]{Theorem}
\newtheorem{prop}[df]{Proposition}

\newtheorem{ex}[df]{Example}

\newtheorem{lem}[df]{Lemma}

\newcommand{\pf}{\noindent{\sc Proof.}\ }

\def\boom{\quad\lower3pt\hbox{\vrule height1.1ex width .9ex depth -.2ex}
                    \vskip9pt}

\let\phi=\varphi
\let\Ga=\Gamma
\let\Om=\Omega

\let\da=\partial

\def\CDO{\mathop{\rm CDO}}

\def\id{{\rm id}}

\def\LAgpd{$\mathcal{LA}$-groupoid}

\def\LAvb{$\mathcal{LA}$-vector bundle}

\def\R{\mathbb{R}}
\def\chigh{{\raise1.5pt\hbox{$\chi$}}}

\def\tila{\widetilde a}
\def\tilb{\widetilde b}

\def\tilzero{\widetilde 0}        
\def\tilq{\widetilde q}
\def\tilx{\widetilde x}

\newcommand{\gog}{\mathfrak{g}}

\let\isom=\cong

\newcommand{\xclam}{^{\textstyle !}}

\let\sdp=\ltimes
\let\pds=\rtimes

\let\sol=\bullet

\def\pback#1{\mathbin{\times_{#1}}}
\def\Whitney#1{\mathbin{\oplus_{#1}}}

\def\llangle{\langle\!\langle}
\def\rrangle{\rangle\!\rangle}

\let\Bar=\overline
\let\Hat=\widehat
\let\Tilde=\widetilde

\def\hatt{\Hat{\phantom{X}}}

\def\dpl{+\hskip-6pt +\hskip4pt}
\def\dminus{\raise2pt\hbox{\vrule height1pt width 2ex}\hskip3pt}
\def\dtimes{\mathbin{\hbox{\huge.}}}

\def\plusH{\ \lower 5pt\hbox{${\buildrel {\textstyle +}
\over {\scriptscriptstyle H}}$}\ }
\def\minusH{\ \lower 5pt\hbox{${\buildrel {\textstyle -}
\over {\scriptscriptstyle H}}$}\ }
\def\timesH{\ \lower 4pt\hbox{${\buildrel {\textstyle .}
\over{\scriptscriptstyle H}}$}\ }

\def\plusV{\ \lower 5pt\hbox{${\buildrel {\textstyle +}
\over {\scriptscriptstyle V}}$}\ }
\def\minusV{\ \lower 5pt\hbox{${\buildrel {\textstyle -}
\over {\scriptscriptstyle V}}$}\ }
\def\timesV{\ \lower 4pt\hbox{${\buildrel {\textstyle .}
\over{\scriptscriptstyle V}}$}\ }

\def\plusK{\ \lower 5pt\hbox{${\buildrel {\textstyle +}
\over {\scriptscriptstyle K}}$}\ }
\def\minusK{\ \lower 5pt\hbox{${\buildrel {\textstyle -}
\over {\scriptscriptstyle K}}$}\ }
\def\timesK{\ \lower 4pt\hbox{${\buildrel {\textstyle .}
\over{\scriptscriptstyle K}}$}\ }

\def\ts#1{\stackrel{\vee}{#1}}

\def\hcompo#1#2#3{{\vcenter{\vbox{\hrule height.#2pt\hbox{\vrule width.#2pt
   height#1pt\kern#3pt\vrule width.#2pt\kern#3pt\vrule width.#2pt}
   \hrule height.#2pt}}}}

\def\vcompo#1#2#3{{\vcenter{\vbox{\hrule height.#2pt
   \hbox{\vrule width.#2pt height#3pt\kern#1pt\vrule width.#2pt}
   \hrule height.#2pt
   \hbox{\vrule width.#2pt height#3pt\kern#1pt\vrule width.#2pt}
   \hrule height.#2pt}}}}

\def\dsq{\mathop{\lower1pt\vbox{\hrule height.4pt \hbox
{\vrule width.4pt height.6em
\kern.6em \vrule width.4pt} \hrule height.4pt}}}

\def\dcomp{\mathop{\dsq\hskip-.88em\raise1pt\hbox{$\scriptstyle\nwarrow$}}}

\def\sq{\vbox{\hrule height.4pt \hbox{\vrule width.4pt height1in
\kern1in \vrule width.4pt} \hrule height.4pt}}

\def\ssq{\vbox{\hrule height.4pt \hbox{\vrule width.4pt height.7in
\kern.7in \vrule width.4pt} \hrule height.4pt}}

\def\tsq{\mathop{\lower1pt\vbox{\hrule height.4pt \hbox
{\vrule width.4pt height.7em
\kern.7em \vrule width.4pt} \hrule height.4pt}}}

\def\sgpd{\,\lower1pt\hbox{$\mlra$}\hskip-0.4in\raise2pt\hbox{$\mlra$}\,}

\def\vgpd{\Bigg\downarrow\!\!\Bigg\downarrow}

\def\upa{\uparrow}

\renewcommand{\to}{\longrightarrow}
\let\mapsto=\longmapsto

\def\longdownmapsto{\raise4ex\hbox{$\Big\downarrow$}}

\def\longdownmaps2{\makebox[0cm]{--}\makebox[0cm]{\raise-2.45ex\hbox{$\Big\downarrow$}}}

\def\gpd{\,\lower1pt\hbox{$\longrightarrow$}\hskip-.24in\raise2pt
             \hbox{$\longrightarrow$}\,}

\def\mlra{\hbox{$\,-\!\!\!-\!\!\!\longrightarrow\,$}}

\def\vgpd{\Bigg\downarrow\!\!\Bigg\downarrow}

\begin{document}

\title{{\bf NOTIONS OF DOUBLE FOR LIE ALGEBROIDS}
\thanks{2000 {\em Mathematics Subject Classification.}
Primary 17B62. Secondary 17B66, 18D05, 22A22, 53D17, 58H05.}
}

\author{K. C. H. Mackenzie\\
        Department of Pure Mathematics\\
        University of Sheffield\\
        Sheffield, S3 7RH\\
        United Kingdom\\
        {\sf K.Mackenzie@sheffield.ac.uk}}

\date{November 5, 2000}

\maketitle

\newpage

\section{INTRODUCTION}

The notion of double introduced by Drinfel'd \cite{Drinfeld:1987} can
be justified by philosophical considerations which lie at the heart of
quantum theory. More concretely, the quantum double of a Hopf algebra
and the classical double of a Lie bialgebra enable the complicated
relations inherent in the original structures to be expressed in simpler
terms on a larger object.

In this paper we are concerned with a specific aspect of the classical
double of a Lie bialgebra: that a pair of Lie algebra structures on a
vector space $\gog$ and its dual $\gog^*$ form a Lie bialgebra if and
only if $\gog\oplus\gog^*$ admits a Lie algebra structure making
$(\gog\oplus\gog^*, \gog, \gog^*)$ a Manin triple. The structure on
$\gog\oplus\gog^*$ may also be viewed as providing a concise formulation
of the twisted derivation equations
\begin{eqnarray*}                                   \label{eq:first}
\rho_X([Y_1, Y_2]) = [\rho_X(Y_1), Y_2] + [Y_1, \rho_X(Y_2)]
    +\rho_{\sigma_{Y_2}(X)}(Y_1) - \rho_{\sigma_{Y_1}(X)}(Y_2),\\
\sigma_Y([X_1, X_2]) = [\sigma_Y(X_1), X_2] + [X_1, \sigma_Y(X_2)]
    +\sigma_{\rho_{X_2}(Y)}(X_1) - \sigma_{\rho_{X_1}(Y)}(X_2),
\end{eqnarray*}
where $X_1, X_2\in{\mathfrak g},\ Y_1, Y_2\in{\mathfrak g}^*$.
Here $\rho$ and $\sigma$ are the two coadjoint representations,
$\rho$ of ${\mathfrak g}$ on ${\mathfrak g}^*$ and $\sigma$ of
${\mathfrak g}^*$ on ${\mathfrak g}$. The equations hold if and
only if $({\mathfrak g}, {\mathfrak g}^*)$, or equivalently
$({\mathfrak g}^*, {\mathfrak g})$, is a Lie bialgebra.

The Manin triple characterization of Lie bialgebras has been extended in
two distinct directions. Firstly, several authors have independently
extended the Manin triple theorem to pairs of Lie algebras not in duality.
The notion of a pair of Lie algebras ${\mathfrak g}, {\mathfrak h}$ such
that the vector space direct sum ${\mathfrak g}\oplus{\mathfrak h}$ has
a Lie algebra structure for which ${\mathfrak g}$ and ${\mathfrak h}$
are Lie subalgebras was introduced in this context by Kosmann--Schwarzbach
and Magri \cite{Kosmann-SchwarzbachM:1988} under the name {\em extension
bicrois\'ee} or {\em twilled extension}, by Lu and Weinstein \cite{LuW:1990}
as {\em double Lie algebra}, and by Majid \cite{Majid:1990}, with the term
{\em matched pair} which we use. (In fact, forms of the concept had been
found much earlier; see \cite{Weinstein:1990}, \cite{ChariP:1994} for
references.) A matched pair structure on a pair of Lie algebras
${\mathfrak g}, {\mathfrak h}$ is equivalent to a pair of representations,
of ${\mathfrak g}$ on ${\mathfrak h}$ and of ${\mathfrak h}$ on
${\mathfrak g}$, which satisfy the twisted derivation equations above.

In the case of Lie bialgebras, the notion of Poisson Lie group is a
natural global concept \cite{Drinfeld:1983}. It was shown by Lu and
Weinstein \cite{LuW:1989} that a further step is possible: every Poisson
Lie group may be integrated in a suitable sense to a double symplectic
groupoid. The appearance of groupoids here may be understood as an
instance of integrability for Poisson manifolds\cite{Karasev:1989},
\cite{Weinstein:1987}, \cite{CosteDW}: the symplectic double groupoid
$S$ of \cite{LuW:1989} corresponding to a Lie bialgebra
$({\mathfrak g}, {\mathfrak g}^*)$ incorporates two (symplectic) groupoid
structures on $S$ which realize the Poisson structures on $G$ and $G^*$,
the Poisson Lie groups which correspond to ${\mathfrak g}$ and ${\mathfrak g}^*$.

Under suitable conditions this symplectic double groupoid is defined by two
global actions, of $G$ on $G^*$ and of $G^*$ on $G$, which integrate the
infinitesimal dressing transformations; this gives $(G, G^*)$ the structure of
a matched pair of Lie groups as defined by \cite{LuW:1990}, \cite{Majid:1990}.
More generally, both \cite{LuW:1990} and \cite{Majid:1990} gave mild conditions
for the global integrability of a matched pair of Lie algebras to a matched pair
of Lie groups, with \cite{Majid:1990} using mainly the actions and \cite{LuW:1990}
mainly the structure of the double.

In fact groupoid structures are present in any matched pair of Lie
groups $(G, H)$: the action of $H$ on $G$ defines a groupoid structure
on $S = G\times H$ with base $G$, and the action of $G$ on $H$ defines a
groupoid structure on $S$ with base $H$. These constitute a double
groupoid structure on $S$, in the categorical sense of Ehresmann; that is,
each groupoid structure on $S$ makes it a
groupoid object in the category of all groupoids. In this case there is
the further property that elements of $S$, regarded as squares, are determined
by any two adjacent sides. It was shown in \cite[\S2]{Mackenzie:1992} that
any double groupoid with this property --- there called {\em vacancy} ---
determines a matched pair structure on its side groupoids. One consequence
of \cite{LuW:1989} is that a matched pair of Lie algebras which is not
globally integrable to a matched pair of Lie groups may nonetheless
integrate to a double groupoid.

The notion of matched pair has been extended to Lie groupoids
\cite[\S2]{Mackenzie:1992} and Lie algebroids \cite{Mokri:1997}; Lu
\cite{Lu:1997} demonstrated that any Poisson Lie group action gives rise
to a matched pair of Lie algebroids which underlies Drinfel'd's
classification of Poisson homogeneous spaces \cite{Drinfeld:1993}.

The second direction in which the Manin triple theorem may be extended is
that of Poisson manifolds which are not groups. A Poisson structure on
a general manifold $P$ gives rise to a Lie algebroid structure on $T^*P$
which, together with the standard tangent bundle structure on $TP$, forms
a Lie bialgebroid in the sense of the author and Ping Xu
\cite{MackenzieX:1994}. In general a Lie bialgebroid $(A, A^*)$ is a Lie
algebroid $A$ together with a Lie algebroid structure on $A^*$ such that
(as \cite{Kosmann-Schwarzbach:1995} showed) the coboundary induced by
each is a derivation of the Schouten bracket of the other; thus Lie
bialgebroids are strong differential Gerstenhaber algebras for which the
underlying module structure comes from a vector bundle \cite{Xu:1999}.
Lie bialgebroids arose in \cite{MackenzieX:1994} as the infinitesimal
invariants of the Poisson groupoids of Weinstein \cite{Weinstein:1988}
(see also \cite{Xu:1995}); they have also arisen from the
Poisson--Nijenhuis structures of Kosmann--Schwarzbach and Magri
\cite{Kosmann-SchwarzbachM:1990} and from the work on dynamical CYBE by
Etingof and Varchenko \cite{EtingofV:1998}. In the latter case
specialization of the algebraic structure leads to quite intricate
systems of pde.

In the case of general Lie bialgebroids there are no representations
corresponding to the coadjoint representations $\rho$ and $\sigma$ above.
One can see immediately that a straightforward extension of the Drinfel'd
double formula to Lie bialgebroids can be made only when the bialgebroid
is in fact a bialgebra.

Recently Liu, Weinstein and Xu \cite{LiuWX:1997} have given an intricate
and highly nonobvious extension of the Manin triple structure on the
double of a Lie bialgebra to a bracket on $\Ga(A\oplus A^*)$, for an
arbitrary Lie bialgebroid $(A, A^*)$. The resulting structure is not usually
a Lie algebroid, but its properties are abstracted by \cite{LiuWX:1997} into
the notion of Courant algebroid; \cite{LiuWX:1997} shows that the Courant
algebroid of a Lie bialgebroid is characterized by the existence of two
complementary integrable isotropic subbundles. Thus the main criterion for
a double is met.

Very recently, Roytenberg \cite{Roytenberg:thesis} has given a treatment
of concepts of double in terms of super geometry.

In this paper we define a very general notion of ``double Lie algebroid''
which includes both the double of a matched pair of Lie algebroids, and
the double of a Lie bialgebroid --- in the bialgebra case, both reduce
to the standard Manin triple theorem. Notice that these two generalizations
of Lie bialgebras do not overlap except in the actual bialgebra case. At the
same time, this notion of double incorporates such objects as the double
(or iterated) tangent bundle $T^2M = T(TM)$ of a smooth manifold $M$,
and other iterated tangent and cotangent structures; these, together with
the canonical isomorphisms between them, are basic to some approaches to
mechanics \cite{Tulczyjew}. The word ``double'' here carries several
meanings: it refers to the classical Drinfel'd double, as derived from
the matched pair structure of a Lie bialgebra, and as derived from the
duality between the Poisson structures in a Lie bialgebra, and it also
refers to the categorical, or iterated, concept of double. The results
of this paper in particular provide global objects corresponding to the
double in all the cases discussed above.

A double Lie algebroid is first of all a double vector bundle as in
Figure \ref{fig:intro}(a);
\begin{figure}[htb]
$$
\begin{matrix}
        &&&&\cr
          &\mathcal{A}&\mlra &A^V&\cr
          &&&&\cr
          &\Bigg\downarrow& &\Bigg\downarrow&   \cr
          &&&&\cr
          &A^H&\mlra &M&\cr
          &&&&\cr
          &&\mbox{(a)}&&\cr
\end{matrix}
\qquad\qquad
\begin{matrix}
                      &&   &&\cr
                            &TA    &\mlra &TM&\cr
                            &&&&\cr
                            &\Bigg\downarrow&&\Bigg\downarrow&\cr
                            &&&&\cr
                            &A &\mlra &M&\cr
                            &&&&\cr
                            &&\mbox{(b)}\cr
\end{matrix}
\qquad\qquad
\begin{matrix}
      &&&&\cr
        &T^*A&\mlra   &A^*&\cr
        &&&&\cr
        &\Bigg\downarrow& &\Bigg\downarrow&   \cr
        &&&&\cr
        &A&\mlra &M&\cr
        &&&&\cr
        &&\mbox{(c)}\cr
\end{matrix}
$$
\caption{\ \label{fig:intro}}
\end{figure}
that is, $\mathcal{A}$ has two vector bundle structures, on bases $A^V$
and $A^H$, each of which is itself a vector bundle on base $M$, such that
for each structure on $\mathcal{A}$, the structure maps (projection,
addition, scalar multiplication) are vector bundle morphisms with respect
to the other structure (see \cite{Pradines:DVB} or
\cite[\S1]{Mackenzie:1992}). Two cases to keep in mind are the tangent
prolongation of an ordinary vector bundle as in Figure~\ref{fig:intro}(b)
(see \cite{Besse} for a classical treatment), and the cotangent double
vector bundle as in Figure~\ref{fig:intro}(c) (see \cite{MackenzieX:1994}).

Now suppose that all four sides of Figure~\ref{fig:intro}(a) have Lie
algebroid structures. The problem is to define compatibility between the
bracket structures on $\mathcal{A}$, bearing in mind that the brackets
are on different modules of sections. The key is the duality for double
vector bundles introduced by Pradines \cite{Pradines:1988}; see
\S\ref{sect:ctvb}. Although we do not know how to interpret our compatibility
condition as a direct condition of the form `either bracket structure on
(sections of) $\mathcal{A}$ is a morphism with respect to the other', we
show in Theorem~\ref{thm:infinv} that this concept of double Lie algebroid
includes the infinitesimal invariants of double Lie groupoids (and the
\LAgpd s of \cite[\S4]{Mackenzie:1992}), which {\em are} defined by such
conditions. The results of \S\ref{sect:adla} were announced in
\cite{Mackenzie:1998}, which also provides an overview of the background
to the present paper.

In \S\ref{sect:dlalba} we extend to abstract Lie bialgebroids the Manin
triple characterization of Lie bialgebras. A Lie algebroid structure on
a vector bundle $A$ induces a Poisson structure on its dual $A^*$ and this
in turn induces a Lie algebroid structure on $T^*A^*\to A^*$. If $A^*$ also
has a Lie algebroid structure, {\em a priori} unrelated, then $T^*A\to A$
likewise has a Lie algebroid structure. Using the canonical isomorphism
$R\colon T^*A^*\to T^*A$ which interchanges $A$ and $A^*$
\cite{MackenzieX:1994}, these equip the cotangent double vector bundle of
Figure~\ref{fig:intro}(c) with four Lie algebroid
structures. We prove in Theorem \ref{theorem:Manin} that these constitute
a double Lie algebroid if and only if $(A, A^*)$ is a Lie bialgebroid.
In the case of Lie bialgebras, \ref{theorem:Manin} is equivalent to the
Manin triple characterization, but in the general case it is necessary to
use the side structures, with bases $A$ and $A^*$, rather than a
structure over $M$. Compared with the general Manin triple theorem of
\cite{LiuWX:1997}, the notion of double Lie algebroid involves a pair of
(relatively) simple bracket structures, the relationship between which
embodies the data, rather than the single but more complicated structure
of a Courant algebroid. The results of \S\ref{sect:dlalba} have been
summarized in \cite{Mackenzie:1998}.

In \S\ref{sect:mpvdla} we show that matched pair structures on a pair
of Lie algebroids $A, B$ with base $M$ correspond to double Lie algebroid
structures on $A\pback{M}B$. Here $A\pback{M}B$ denotes the double vector
bundle formed by the two pullbacks $q_A^!B$ and $q_B^!A$ across the
projections $q_A\colon A\to M$ and $q_B\colon B\to M$; we call such double
vector bundles {\em vacant}. The situation in \S\ref{sect:mpvdla} is that
in addition to the Lie algebroid structures on $A$ and $B$ over $M$, we
have Lie algebroid structures on the two pullbacks $q_A^!B$ and $q_B^!A$,
subject to the conditions of \S\ref{sect:adla}. It does not seem evident
that such a pair of structures is equivalent to a pair of representations
satisfying the twisted derivation conditions (together with a third which
is vacuous when $M$ is a point; see \ref{df:mokri}), and thus in turn
equivalent to a single Lie algebroid structure on the Whitney sum $A\oplus B$
as in \ref{prop:mokri}. We prove this in \ref{theorem:mp}.

It then follows in Theorem \ref{theorem:sdp} that a matched pair structure
on Lie algebroids $A, B$ with given representations $\rho, \sigma$ is
equivalent to a Lie bialgebroid structure on $(A\sdp B^*, A^*\pds B)$, where
the semi--direct products are formed with respect to the contragredient
representations.
Thus the single Lie bialgebroid equation, applied to these
semi--direct products, also specializes to the system of matched pair
equations.

One immediately striking feature of the notion of Lie bialgebra is the
number of different perspectives from which it can be understood
--- as a pair of compatible actions, as a Manin triple, as the
infinitesimal structure of a Poisson group.
We show in this paper that these different approaches, extended to
general Lie algebroids and encompassing a very much broader range of
phenomena, can be unified into the single concept of double Lie algebroid.
Because the bialgebroid equation is the key compatibility condition for
a double Lie algebroid, and may be formulated diagrammatically, all
of the notions subsumed under that of double Lie algebroid should be capable
of formulation in fairly general categories, including those arising
in modern forms of quantization.

Our definition of double Lie algebroid depends crucially on the
underlying duality for finite rank vector bundles. It would be very
interesting to be able to define a comparable concept of double for
Lie algebroid--like structures on general modules (variously known
as Lie pseudoalgebras, Gerstenhaber algebras, Lie--Rinehart algebras,
and a number of other terminologies).

A subsequent article \cite{Mackenzie:UA} extends to actions of Poisson
groupoids Lu's construction \cite{Lu:1997} of a matched pair of Lie
algebroids for a Poisson group action. In the general case of Poisson
groupoid actions, this produces a double Lie algebroid which does not
correspond to either a matched pair or a Lie bialgebroid.

I am very grateful to Yvette Kos\-mann--Schwarz\-bach, Alan Weinstein and
Ping Xu for conversations as this material developed over a number of years.

\section{DOUBLE VECTOR BUNDLES AND THEIR DUALITY}
\label{sect:ctvb}

We begin by recalling the duality of double vector bundles from
\cite{Mackenzie:1999} or \cite{KoniecznaU:1999}. We then show that the
existence of the duality between the duals of a double vector bundle
reflects the existence of triple structures associated with its tangent
and cotangent. The use of these triple structures shortens many arguments
throughout the paper.

The definition of {\em double vector bundle} given in the Introduction (see
Figure~\ref{fig:intro}(a)) is complete but needs elaboration. The
commutativity conditions stated there are precisely what is needed to ensure
that when four elements $\xi_1,\dots,\xi_4$ of $\mathcal{A}$ are such that
the LHS of
$$
(\xi_1\plusH\xi_2)\plusV(\xi_3\plusH\xi_4) =
(\xi_1\plusV\xi_3)\plusH(\xi_2\plusH\xi_4)
$$
is defined, then the RHS is also, and they are equal. Here $\plusH$ denotes
the addition in $\mathcal{A}\to A^V$ and $\plusV$ the addition in
$\mathcal{A}\to A^H$. See \cite{Pradines:DVB} or \cite[\S1]{Mackenzie:1992}.

We denote the projections in Figure~\ref{fig:intro}(a) by
$\tilq_H\colon\mathcal{A}\to A^V$ and $\tilq_V\colon\mathcal{A}\to A^H$;
these are vector bundle morphisms over the projections $q_H\colon A^H\to M$
and $q_V\colon A^V\to M$ respectively. The intersection of the kernels of
the two projections defined on $\mathcal{A}$ is called the {\em core} of
the double vector bundle $\mathcal{A}$ and denoted $K$. The vector bundle
structures on $\mathcal{A}$ induce a common vector bundle structure on $K$
with base $M$. The kernel of $\tilq_H\colon\mathcal{A}\to A^V$ is now the
pullback $q_H\xclam K$ of $K$ across $q_H$, and the kernel of
$\tilq_V\colon\mathcal{A}\to A^H$ is $q_V\xclam K$. In practice cores are
often identified with structures which are not strictly subsets of
$\mathcal{A}$ and in such cases we write $\Bar{k}$ for the element of
$\mathcal{A}$ corresponding to a $k\in K$.

Given $X\in A^H$ we denote the zero of $\mathcal{A}\to A^H$ above $X$
by $\tilzero^V_X$. Similarly the zero of $\mathcal{A}\to A^V$ above
$x\in A^V$ is denoted $\tilzero^H_x$.

Consider the dual $\mathcal{A}^{*V}$ of the vertical bundle structure
$\mathcal{A}\to A^H$. In addition to its standard structure on base $A^H$,
this has a vector bundle structure on base $K^*$, the projection of which
is defined by
\begin{equation}                             \label{eq:unfproj}
\langle\tilq^{(*)}_V(\Phi), k\rangle =
      \langle\Phi, \tilzero^V_X \plusH \Bar k\rangle
\end{equation}
where $\Phi\colon \tilq_V^{-1}(X)\to\R,\ X\in A_m^H,$ and $k\in K_m$.
The addition in $\mathcal{A}^{*V}\to K^*$, which we also denote by $\plusH$,
is defined by
\begin{equation}                              
\langle\Phi\plusH\Phi', \xi\plusH\xi'\rangle =
   \langle\Phi,\xi\rangle + \langle\Phi',\xi'\rangle.
\end{equation}
The zero of $\mathcal{A}^{*V}\to K^*$ above $\kappa\in K^*_m$ is
denoted $\tilzero^{(*V)}_\kappa$ and defined by
$$
\langle\,\tilzero^{(*V)}_\kappa, \tilzero^H_x \plusV \Bar k\rangle =
\langle\kappa, k\rangle
$$
where $x\in A^V_m, k\in K_m.$ The scalar multiplication is defined in a
similar way. These two structures make $\mathcal{A}^{*V}$ into a double
vector bundle as in Figure~\ref{fig:dvbduals}(a),
\begin{figure}[hbt]
$$
\begin{matrix}
      && \tilq^{(*)}_V &&\cr
        &\mathcal{A}^{*V}&\mlra   &K^*&\cr
        &&&&\cr
\tilq_{*V} &\Bigg\downarrow& &\Bigg\downarrow&q_{K^*} \cr
        &&&&\cr
        &A^H&\mlra &M&\cr
        &&q_H&&\cr
        &&&&\cr
        &&\mbox{(a)}&&\cr
\end{matrix}
\qquad\qquad\qquad
\begin{matrix}
      && \tilq_{*H}  &&\cr
        &\mathcal{A}^{*H}&\mlra   &A^V&\cr
        &&&&\cr
\tilq^{(*)}_H &\Bigg\downarrow& &\Bigg\downarrow&q_V \cr
        &&&&\cr
        &K^*&\mlra &M&\cr
        &&q_{K^*}&&\cr
        &&&&\cr
        &&\mbox{(b)}&&\cr
\end{matrix}
$$
\caption{\ \label{fig:dvbduals}}
\end{figure}
the {\em vertical dual of $\mathcal{A}$}. The core of $\mathcal{A}^{*V}$
identifies with $(A^V)^*$, with the core element $\Bar\psi$
corresponding to $\psi\in(A^V_m)^*$ given by
$$
\langle\Bar\psi, \tilzero^H_x \plusV \Bar k\rangle =
\langle\psi,x\rangle.
$$
See \cite{Pradines:1988} or \cite[\S3]{Mackenzie:1999}.

There is also a {\em horizontal dual} $\mathcal{A}^{*H}$ with sides $A^V$
and $K^*$ and core $(A^H)^*$, as in Figure~\ref{fig:dvbduals}(b). There is
now the following somewhat surprising result.

\begin{thm}[{\cite[3.1]{Mackenzie:1999}}, {\cite[Thm.~16]{KoniecznaU:1999}}]
\label{thm:dualduality}
There is a natural (up to sign) duality between $\mathcal{A}^{*V}\to K^*$ and
$\mathcal{A}^{*H}\to K^*$ given by
\begin{equation}                       \label{eq:3duals}
\langle\Phi, \Psi\rangle = \langle\Psi, \xi\rangle
                            -  \langle\Phi, \xi\rangle
\end{equation}
where $\Phi\in \mathcal{A}^{*V},\ \Psi\in \mathcal{A}^{*H}$ have
$\tilq_V^{(*)}(\Phi) = \tilq_H^{(*)}(\Psi)$ and $\xi$ is any element
of $\mathcal{A}$ with $\tilq_V(\xi) = \tilq_{*V}(\Phi)$ and
$\tilq_H(\xi) = \tilq_{*H}(\Psi).$
\end{thm}

The pairing on the LHS of (\ref{eq:3duals}) is over $K^*$, whereas the
pairings on the RHS are over $A^V$ and $A^H$ respectively.

\begin{ex}\rm
Given an ordinary vector bundle $(A,q,M)$ there is the tangent double
vector bundle of Figure~\ref{fig:intro}(b); for a detailed account of this
see \cite{Besse} or \cite[\S5]{MackenzieX:1994}. The core of $TA$ consists
of the vertical vectors along the zero section and identifies canonically
with $A$. Dualizing the standard structure $TA\to A$ gives the cotangent
double vector bundle of Figure~\ref{fig:intro}(c) with core $T^*M$. Dualizing
the structure $TA\to TM$ gives a double vector bundle which we denote
$(T^\sol A;A^*,TM;M)$; this is canonically isomorphic to $(T(A^*);A^*,TM;M)$
under an isomorphism $I\colon T(A^*)\to T^\sol A$ given by
$$
\langle I(\mathcal{X}),\xi \rangle =
\llangle\mathcal{X}, \xi \rrangle
$$
where $\mathcal{X}\in T(A^*)$ and $\xi \in TA$ have
$T(q_*)(\mathcal{X}) = T(q)(\xi)$, and $\llangle\ ,\ \rrangle$ is the
tangent pairing of $T(A^*)$ and $TA$ over $TM$.
See \cite[\S5]{MackenzieX:1994}.

The structure of $\mathcal{A} = TA$ thus induces a pairing of $T^*A$ and
$T(A^*)$ over $A^*$ given by
$$
\langle\Phi,\mathcal{X}\rangle =
\llangle\mathcal{X},\xi\rrangle -
\langle\Phi, \xi\rangle
$$
where $\Phi\in T^*A$ and $\mathcal{X}\in T(A^*)$ have
$r(\Phi) = p_{A^*}(\mathcal{X})$, and $\xi\in TA$ is chosen so that
$T(q)(\xi) = T(q_*)(\mathcal{X})$ and $p_A(\xi) = c_A(\Phi)$.
(Here $c_A$ is the projection of the standard cotangent bundle and
$r\colon T^*A\to A^*$ is a particular case of (\ref{eq:unfproj}).)
This pairing is nondegenerate by a general result \cite[3.1]{Mackenzie:1999}
so it defines an isomorphism of double vector bundles
$R\colon T^*A^*\to T^*A$ by the condition
$$
\langle R(\mathcal{F}), \mathcal{X}\rangle =
            \langle\mathcal{F}, \mathcal{X}\rangle
$$
where the pairing on the RHS is the standard one of $T^*(A^*)$ and
$T(A^*)$ over $A^*$. This $R$ preserves the side bundles $A$ and $A^*$
but induces $-\id\colon T^*M\to T^*M$ as the map of cores. It is an
antisymplectomorphism with respect to the exact symplectic structures.
In summary we now have the very useful equation
\begin{equation}                                   \label{eq:vue}
\langle\mathcal{F}, \mathcal{X}\rangle + \langle R(\mathcal{F}), \xi\rangle
= \llangle\mathcal{X}, \xi\rrangle,
\end{equation}
for $\mathcal{F}\in T^*A^*,\ \mathcal{X}\in T(A^*),\ \xi\in TA$, where the
pairings are over $A^*, A$ and $TM$ respectively. For all of this see
\cite[\S5]{MackenzieX:1994}.
\end{ex}

Now return to the general double vector bundle in Figure~\ref{fig:intro}(a),
denoting the core by $K$. Each vector bundle structure on $\mathcal{A}$
gives rise to a double cotangent bundle. These two double cotangents
fit together into a triple structure as in Figure~\ref{fig:cottrip}(a).
This structure is a triple vector bundle in an obvious sense: each face
is a double vector bundle and for each vector bundle structure on
$T^*\mathcal{A}$, the maps defining the structure are morphisms of double
vector bundles. The left and rear faces of Figure~\ref{fig:cottrip}(a)
are the two cotangent doubles of the two structures on $\mathcal{A}$ and
the top face may be regarded as the cotangent double of either of the two
duals of $\mathcal{A}$. (In all diagrams of this type, we take the oblique
arrows to be coming out of the page.) The fact that the cube in
Figure~\ref{fig:cottrip}(a) may be rotated about the diagonal from
$T^*\mathcal{A}$ to $M$, preserving its type, is essentially a consequence
of Theorem~\ref{thm:dualduality}.

\begin{figure}[htb]
\begin{picture}(340,180)(-20,30)
\put(-10,150){$\begin{matrix}
                &&      &\cr
                      &T^*\mathcal{A}&\mlra &\mathcal{A}^{*H}\cr
                      &&&\cr
                      &\Bigg\downarrow   & &\Bigg\downarrow     \cr
                      &&&\cr
                      &\mathcal{A}&\mlra &A^V\cr
                \end{matrix}$}

\put(25, 165){\vector(1,-1){20}}                
\put(95, 165){\vector(1,-1){20}}               
\put(25, 100){\vector(1,-1){20}}                
\put(95, 100){\vector(1,-1){20}}               

\put(45,88){$\begin{matrix}
               &&      &\cr
                     &\mathcal{A}^{*V} &\mlra &K^*\cr
                     &&&\cr
                     &\Bigg\downarrow &&\Bigg\downarrow \cr
                     &&&\cr
                     &A^H &\mlra & M\cr
                     &&&\cr
                     &&\mbox{(a)}&\cr
               \end{matrix}$}


\put(220,150){$\begin{matrix}
                &&      &\cr
                      &T\mathcal{A}&\mlra &T(A^V)\cr
                      &&&\cr
                      &\Bigg\downarrow & &\Bigg\downarrow \cr
                      &&&\cr
                      &\mathcal{A}&\mlra &A^V \cr
                \end{matrix}$}

\put(255, 165){\vector(1,-1){20}}                
\put(325, 165){\vector(1,-1){20}}               
\put(255, 100){\vector(1,-1){20}}                
\put(325, 100){\vector(1,-1){20}}               

\put(265,88){$\begin{matrix}
               &&      &\cr
                     &T(A^H) &\mlra &TM\cr
                     &&&\cr
                     &\Bigg\downarrow &&\Bigg\downarrow \cr
                     &&&\cr
                     &A^H &\mlra & M\cr
                     &&&\cr
                     &&\mbox{(b)}&\cr
               \end{matrix}$}
\end{picture}\caption{\ \label{fig:cottrip}}
\end{figure}

Figure~\ref{fig:cottrip}(a) is, in a sense we will make precise elsewhere,
the vertical dual of the tangent prolongation of $\mathcal{A}$ as given in
Figure~\ref{fig:cottrip}(b). Five of the six faces in \ref{fig:cottrip}(a)
are double vector bundles of types considered already; it is only necessary
to verify that the top face is a double vector bundle. The cores of the same
five faces are known, and we take the core of the top face to be $T^*K$.
(This is a special case of \cite[1.5]{Mackenzie:1999}.)
Taking these cores in pairs, together with a parallel edge, then gives three
double vector bundles: the left--right cores form $(T^*(A^H);A^H,(A^H)^*;M)$,
the back--front cores form $(T^*(A^V);A^V,(A^V)^*;M)$ and the up--down
cores form $(T^*K;K,K^*;M)$. Each of these {\em core double vector bundles}
has core $T^*M$.

We will also need to consider the cotangent triples of the two duals of
$\mathcal{A}$. Figure~\ref{fig:cofd}(a) is the cotangent triple of the double
vector bundle $(\mathcal{A}^{*V};A^H,K^*;M)$; the ${}^\dagger$ denotes the
dual over $K^*$. We use the isomorphisms of double vector bundles
$$
Z_V\colon (\mathcal{A}^{*H})^\dagger \to \mathcal{A}^{*V}, \qquad
Z_H\colon (\mathcal{A}^{*V})^\dagger \to \mathcal{A}^{*H}
$$
induced by the pairing (\ref{eq:3duals}). Note that $Z_V$ preserves both
sides, $A^H$ and $K^*$, but induces $-\id\colon (A^V)^*\to (A^V)^*$ on the
cores, while $Z_H$ preserves $K^*$ and the core $(A^H)^*$, but induces
$-\id$ on the sides $A^V$; this reflects the fact that $Z_V = Z_H^\dagger$,
the dual over $K^*$; see \cite[3.6]{Mackenzie:1999}.

\begin{figure}[h] 
\begin{picture}(340,180)(0,30)               
\put(-10,150){$\begin{matrix}
               &&     &&\cr
                     &T^*(\mathcal{A}^{*V})&\mlra &(\mathcal{A}^{*V})^\dagger & \isom \mathcal{A}^{*H}\cr
                     &&&&\cr
                     &\Bigg\downarrow   & &\Bigg\downarrow  &   \cr
                     &&&&\cr
                     &\mathcal{A}^{*V}&\mlra &K^*&\cr
               \end{matrix}$}

\put(30, 165){\vector(1,-1){20}}                
\put(110, 165){\vector(1,-1){18}}               
\put(30, 100){\vector(1,-1){20}}                
\put(110, 100){\vector(1,-1){20}}               

\put(33,105){$\begin{matrix}
               &&      &\cr
                     &\phantom{X}\mathcal{A}\phantom{X} &\mlra &\phantom{X}A^V\phantom{X}\cr
                     &&&\cr
                     &\Bigg\downarrow &&\Bigg\downarrow \cr
                     &&&\cr
                     &A^H &\mlra & M\cr
               \end{matrix}$}

\put(90,20){(a)}


\put(210,150){$\begin{matrix}
                &&      &\cr
                      &T(\mathcal{A}^{*V})&\mlra &T(K^*)\cr
                      &&&\cr
                      &\Bigg\downarrow & &\Bigg\downarrow \cr
                      &&&\cr
                      &\mathcal{A}^{*V}&\mlra &K^* \cr
                \end{matrix}$}

\put(245, 165){\vector(1,-1){20}}                
\put(325, 165){\vector(1,-1){20}}               
\put(245, 100){\vector(1,-1){20}}                
\put(325, 100){\vector(1,-1){20}}               

\put(255,105){$\begin{matrix}
               &&      &\cr
                     &T(A^H)&\mlra &TM\cr
                     &&&\cr
                     &\Bigg\downarrow &&\Bigg\downarrow \cr
                     &&&\cr
                     &A^H &\mlra & M\cr
               \end{matrix}$}
\put(310,20){(b)}
\end{picture}\caption{\ \label{fig:cofd}}
\end{figure}

When $\mathcal{A} = T^*A$ as in Figure~\ref{fig:intro}(b), we have
$Z_V = R\circ I^\dagger$, where the ${}^\dagger$ dual here is over $A^*$.

\section{PRELIMINARY CASE}
\label{sect:pc}

Before addressing the general concept of double Lie algebroid, it will be
useful to deal with a very special case.

\begin{df}
An {\em \LAvb}\ is a double vector bundle as in Figure~{\rm \ref{fig:intro}(a)}
together with Lie algebroid structures on a pair of parallel sides, such
that the structure maps of the other pair of vector bundle structures are
Lie algebroid morphisms.
\end{df}

For definiteness, take the Lie algebroid structures to be on
$\mathcal{A}\to A^H$ and $A^V\to M$.

In the terminology of \cite[\S4]{Mackenzie:1992}, an \LAvb\ is an \LAgpd\
in which the groupoid structures are vector bundles (and in which the
scalar multiplication also preserves the Lie algebroid structures).
The core of a general \LAgpd\ has a Lie algebroid structure induced from the
Lie algebroid structure on $\mathcal{A}$ \cite[\S5]{Mackenzie:1992}. Namely,
each $k\in\Ga K$ induces $\Bar{k}\in\Ga_{A^H}\mathcal{A}$ defined by
$\Bar{k}(X) = k(m)\plusH\tilzero^V_X$ for $X\in A^H_m$ and the bracket on
$\Ga K$ is obtained by $\Bar{[k,\ell]} = [\Bar{k}, \Bar{\ell}]$.

\begin{lem}                                    \label{lem:K0}
The anchor and the bracket on the core of an \LAvb\ are zero.
\end{lem}

\pf
The anchor $\tila_V\colon\mathcal{A}\to T(A^H)$ is a morphism of double vector
bundles and therefore induces a core map $\da_H\colon K\to A^H$ which, by
\cite[\S5]{Mackenzie:1992}, is a Lie algebroid morphism. Since $A^H$ is
abelian, we have $a_K = a_H\circ\da_H = 0$.

Horizontal scalar multiplication by $t\neq 0$ defines a morphism
$\mathcal{A}\to\mathcal{A}$ over $A^H\to A^H$ and therefore induces a map of
sections $t_H\colon \Ga_{A^H}\mathcal{A} \to \Ga_{A^H}\mathcal{A}$; the Lie
algebroid condition then ensures that
$[t_H(\xi), t_H(\eta)] = t_H([\xi,\eta])$
for all $\xi, \eta\in\Ga_{A^H}\mathcal{A}$. Now for $k\in\Ga K$,
$\Bar{tk} = t_H(\Bar{k})$ and so
$\Bar{t[k,\ell]} = [t_H(\Bar{k}), t_H(\Bar{\ell})]
= [\Bar{tk}, \Bar{t\ell}] = \Bar{t^2[k,\ell]}$. Therefore the bracket
must be zero.
\boom

Since $\mathcal{A}\to A^H$ is a Lie algebroid, the dual $\mathcal{A}^{*V}$ has
its dual Poisson structure, which is linear with respect to the bundle
structure over $A^H$. The remainder of this section gives the proof of
the following result.

\begin{thm}                                     \label{theorem:bothduals}
The Poisson structure on $\mathcal{A}^{*V}$ is also linear with respect to
the bundle structure over $K^*$.
\end{thm}

This is actually a special case of \cite[3.14]{Mackenzie:1999} (providing
attention is paid to the scalar multiplication). We give a direct proof
however, since there are special features which are needed later.

The functions $\mathcal{A}^{*V}\to\R$ which are linear with respect to $K^*$
are determined by sections of $\mathcal{A}^{*H}\to K^*$ via the duality
(\ref{eq:3duals}) which we here denote by $\langle\ ,\ \rangle_{K^*}$
for clarity. Given $\Psi\in\Ga_{K^*}\mathcal{A}^{*H}$ define
$$
\ell^\dagger_\Psi\colon\mathcal{A}^{*V}\to\R, \qquad
\Phi\mapsto\langle\Phi, \Psi(\tilq^{(*)}_V(\Phi))\rangle_{K^*}.
$$
There are two principal cases.
Firstly, consider sections $\xi\in\Ga_{A^H}\mathcal{A},\ x\in\Ga A^V$. The
pair $(\xi, x)$ is a {\em (vertical) linear section} if $(\xi, x)$ is a
vector bundle morphism from $(A^H, q_H, M)$ to $(\mathcal{A}, \tilq_H, A^V)$.

\begin{prop}
{\rm (i)} If $(\xi, x)$ is a linear section, then $\ell_\xi\colon
\mathcal{A}^{*V}\to\R$ is linear with respect to $K^*$ as well as $A^H$, and
the restriction of $\ell_\xi$ to the core of $\mathcal{A}^{*V}$ is
$\ell_x\colon (A^V)^*\to\R$.

{\rm (ii)} If the function $\ell_\xi$ defined by some
$\xi\in\Ga_{A^H}\mathcal{A}$ is linear with respect to $K^*$ as well as
$A^H$, then the restriction of $\ell_\xi$ to the core defines a section
$x\in\Ga A^V$, and $(\xi, x)$ is a linear section.
\end{prop}

\pf
(i) Take elements $(\Phi_1; X_1, \kappa; m)$ and $(\Phi_2; X_2, \kappa; m)$
of $\mathcal{A}^{*V}$. Their horizontal sum is of the form
$(\Phi_1\plusH\Phi_2; X_1 + X_2, \kappa; m)$ and so, using first the
linearity of $\xi$ and then the definition of $\plusH$,
\begin{eqnarray*}
\ell_\xi(\Phi_1\plusH\Phi_2)
& = & \langle \Phi_1\plusH\Phi_2, \xi(X_1 + X_2)\rangle =
\langle \Phi_1\plusH\Phi_2, \xi(X_1) \plusH \xi(X_2)\rangle\\
& = & \langle \Phi_1, \xi(X_1)\rangle + \langle \Phi_2, \xi(X_2)\rangle =
\ell_\xi(\Phi_1) + \ell_\xi(\Phi_2).
\end{eqnarray*}
The remainder of (i) is similar.

(ii) Since the core embedding $(A^V)^*\to\mathcal{A}^{*V}$ is linear with
respect to either structure on $\mathcal{A}^{*V}$, the first statement is
immediate. The relationship between $\xi$ and $x$ is
$\langle\xi(0^H_m), \Bar{\psi}\rangle = \langle x(m), \psi\rangle$
for all $\psi\in(A^V_m)^*$. Writing $y = \tilq_H(\xi(0^H_m))$ we
therefore have $\xi(0^H_m) = \tilzero^H_y\plusV\Bar{k}$ for some
$k\in K_m$. Pairing with any $\Bar{\psi}$ we find $y = x(m)$.

Linearity over $K^*$ means that, for any
$(\Phi_1;X_1, \kappa; m), (\Phi_2;X_2, \kappa; m)$ in $\mathcal{A}^{*V}$,
\begin{equation}                                    \label{eq:linK}
\langle\xi(X_1+X_2), \Phi_1\plusH\Phi_2\rangle =
\langle\xi(X_1), \Phi_1\rangle + \langle\xi(X_2), \Phi_2\rangle.
\end{equation}
Putting $X_1 = X_2 = 0^H_m$ and $\Phi_1 = \Phi_2 = \tilzero_\kappa^{(*V)}$
for any $\kappa\in K^*_m$, we have
$\langle\xi(0^H_m), \tilzero^{(*V)}_\kappa\rangle = 0$. Using
$\xi(0^H_m) = \tilzero_{x(m)}\plusV\Bar{k}$, this gives
$\langle\kappa, k\rangle = 0$ for all $\kappa$, so $k = 0$.

Next in (\ref{eq:linK}) put $\Phi_1 = \Bar{\psi},\ \Phi_2 = \tilzero_X$ and
$X_1 = 0^H_m,\ X_2 = X$. Then we have
\begin{equation*}
\begin{split}
\langle\psi, x(m)\rangle
& = \langle\xi(X), \Bar{\psi}\plusH\tilzero_X\rangle
   = \langle\tilzero^H_y\plusH\xi(X), \Bar{\psi}\plusH\tilzero_X\rangle\\
& = \langle\tilzero^H_y, \Bar{\psi}\rangle + \langle\xi(X), \tilzero_X\rangle
    = \langle\psi, y\rangle + 0
\end{split}
\end{equation*}
where $y = \tilq_H(\xi(X))$. So
$x(q_H(X)) = x(m) = y = \tilq_H(\xi(X))$. The proof that $(\xi, x)$ is
linear is now immediate.
\boom

Given a linear section $(\xi, x)$, denote by $\xi^\sqcap$ the section
of $\mathcal{A}^{*H}\to K^*$ which $\ell_\xi$ defines. Thus
$$
\langle\Phi, \xi^\sqcap(\kappa)\rangle_{K^*} = \ell_\xi(\Phi)
= \langle\Phi, \xi(X)\rangle_{A^H}
$$
for $(\Phi;X, \kappa; m)$ in $\mathcal{A}^{*V}$. The proof of the
result which follows is not difficult.

\begin{prop}
$(\xi^\sqcap, x)$ is a (vertical) linear section of $\mathcal{A}^{*H}$.
\end{prop}

There is thus a bijective correspondence between vertical linear sections
of $\mathcal{A}$ and vertical linear sections of $\mathcal{A}^{*H}$. This of
course applies to any double vector bundle: in the case where $\mathcal{A}$ is
the tangent of an ordinary vector bundle $A$, as in Figure~\ref{fig:intro}(b),
this is the correspondence between linear vector fields on $A$ and linear
vector fields on $A^*$ (see \cite[\S2]{MackenzieX:1998}).

Secondly, any $\phi\in\Ga(A^H)^*$ induces a core section $\Bar{\phi}$
of $\mathcal{A}^{*H}\to K^*$ by
$$
\Bar{\phi}(\kappa) = \tilzero^{(*H)}_\kappa \plusH \Bar{\phi(m)},\qquad
\kappa\in K^*_m,
$$
which induces $\ell^\dagger_{\Bar{\phi}}\colon\mathcal{A}^{*V}\to\R$. A
simple calculation shows that
$$
\ell^\dagger_{\Bar{\phi}} = \ell_\phi\circ\tilq_{*V}.
$$

We also need the pullbacks across $\mathcal{A}^{*V}\to K^*$ of functions
on $K^*$. In particular, for $k\in\Ga K$ and $f\in C(M)$, we have
$$
\ell_k\circ\tilq^{(*)}_V = \ell_{\Bar{k}},\qquad
f\circ q_{K^*}\circ\tilq^{(*)}_V = f\circ q_H\circ\tilq_{*V}
$$
where $\Bar{k}\in\Ga_{A^H}\mathcal{A}$ is the core section for $k$.

We now have four classes of functions $\mathcal{A}^{*V}\to\R$. The functions
linear over $K^*$ are generated by the $\ell_\xi$ for $(\xi, x)$ linear,
together with the $\ell_\phi\circ\tilq_{*V}$ for $\phi\in\Ga(A^H)^*$. The
pullbacks from $K^*$ are generated by the $\ell_{\Bar{k}},\ k\in\Ga K$,
together with the $f\circ q_H\circ\tilq_{*V}$ for $f\in C(M)$. In the
proof of \ref{theorem:bothduals} we therefore need only consider functions
of these types; we refer to them as types $L_1, L_2, P_1, P_2$
respectively.

\begin{lem}
If $(\xi, x)$ and $(\eta, y)$ are linear sections, then
$([\xi, \eta], [x,y])$ is also.
\end{lem}

\pf
For a section $\xi$ which projects under $\tilq_H$ to a section $x$,
define a section $\xi\oplus\xi$ of $\mathcal{A}\Whitney{A^V}\mathcal{A}\to
A^H\Whitney{M}A^H$ by $(\xi\oplus\xi)(X\oplus Y) = \xi(X)\oplus\xi(Y)$.
Then $(\xi, x)$ is linear if and only if $\xi\oplus\xi$ projects to
$\xi$ under the horizontal addition
$\plusH\colon \mathcal{A}\Whitney{A^V}\mathcal{A}\to \mathcal{A}$,
which is a Lie algebroid morphism over $+\colon A^H\Whitney{M}A^H\to A^H$
by hypothesis. Since the two components of a $\xi\oplus\xi$ depend on each
variable separately, we have
$[\xi\oplus\xi, \eta\oplus\eta] = [\xi, \eta]\oplus[\xi, \eta]$ and the
result follows.
\boom

From $\{\ell_\xi, \ell_\eta\} = \ell_{[\xi, \eta]}$ it follows that
the bracket of two type $L_1$ functions is $L_1$. The bracket of two $L_2$
functions is zero, since they are pullbacks from $A^H$. For the mixed
case,
$$
\{\ell_\xi, \ell_\phi\circ\tilq_{*V}\} =
\tila_V(\xi)(\ell_\phi)\circ\tilq_{*V}.
$$
Now $\tila_V(\xi)$ is a linear vector field on $A^H$ so,
by \cite[(5)]{MackenzieX:1998}, we can define
$$
\tila_V(\xi)(\ell_\phi) = \ell_{D^{(*)}_\xi(\phi)}
$$
where $D^{(*)}_\xi\colon\Ga(A^H)^*\to\Ga(A^H)^*$. In particular
the bracket of an $L_1$ and an $L_2$ is an $L_2$.

Here $D^{(*)}$ is a covariant differential operator in the sense of
\cite[III\S2]{Mackenzie:LGLADG}; that is, $D = D^{(*)}$ is a linear
differential operator of order $\leq 1$ and there is a vector field
$X = a_V(\xi)$ such that $D(f\phi) = fD(\phi) + X(f)\phi$ for all
$f\in C(M)$ and section $\phi$. For any vector bundle $E$ there is a
vector bundle $\CDO(E)$ the sections of which are the covariant
differential operators, and with anchor $D\mapsto X$ and the usual
bracket, $\CDO(E)$ is a Lie algebroid.

For an $L_1$ and a $P_1$ we have $\{\ell_\xi, \ell_{\Bar{k}}\} =
\ell_{[\xi, \Bar{k}]}$. Now $\xi\sim x$ and $\Bar{k}\sim 0$
under $\tilq_H$, so $[\xi, \Bar{k}]\sim 0$. It is therefore $\Bar{k'}$
for some $k'\in\Ga K$ which we denote $C_\xi(k)$; it is easily
checked that $C_\xi\colon\Ga K\to\Ga K$ is again a covariant
differential operator. Thus we have another $P_1$.

For an $L_1$ and a $P_2$ we have
$$
\{\ell_\xi, f\circ q_H\circ\tilq_{*V}\} =
\tila_V(\xi)(f\circ q_H)\circ\tilq_{*V} =
x(f)\circ q_H\circ\tilq_{*V},
$$
again a $P_2$, and for an $L_2$ and a $P_1$,
$$
\{\ell_\phi\circ\tilq_{*V}, \ell_{\Bar{k}}\} =
   -\tila_V(\Bar{k})(\ell_\phi)\circ\tilq_{*V} =
   -\da_H(k)^\upa(\ell_\phi)\circ\tilq_{*V} =
   -\langle\da_H(k), \phi\rangle\circ q_H\circ\tilq_{*V},
$$
a $P_2$. Here $z^\upa$, for a vector field $z$, is the vertical lift of
$z$ as in \cite[\S2]{MackenzieX:1998}. The remaining four cases give zero,
using \ref{lem:K0}.

This completes the proof of Theorem \ref{theorem:bothduals}.
\boom

Since the Poisson structure is linear over $A^H$,
$\pi^{\#V}\colon T^*(\mathcal{A}^{*V})\to T(\mathcal{A}^{*V})$ is a morphism
of double vector bundles for the left faces of Figure~\ref{fig:cofd},
with the corner map $\mathcal{A}\to T(A^H)$ being $\tila_V$ and core map
$-\tila_V^*\colon T^*(A^H)\to \mathcal{A}^{*V}.$ From \ref{theorem:bothduals}
it now follows that $\pi^{\#V}$ is also a morphism of double vector
bundles for the rear faces of Figure~\ref{fig:cofd}. Denote the corner
map $(\mathcal{A}^{*V})^\dagger\to T(K^*)$ by $\chigh_V$; since $\pi^{\#V}$
is skew--symmetric, the core map of the rear faces is
$-\chigh_V^*\colon T^*K^*\to \mathcal{A}^{*V}$. It then follows by a simple
argument (as in \cite[2.3]{Mackenzie:1999}) that $\pi^{\#V}$ is a
morphism of triple vector bundles.

\section{ABSTRACT DOUBLE LIE ALGEBROIDS}
\label{sect:adla}

We now turn to the general notion of double Lie algebroid. Again consider a
double vector bundle as in Figure~\ref{fig:intro}(a). We now assume that
there are Lie algebroid structures on all four sides. The definition
comprises three conditions.

\subsection{Condition I}

{\em With respect to the two vertical Lie algebroids, $\mathcal{A}\to A^H$
and $A^V\to M$, the double vector bundle $\mathcal{A}$ is an \LAvb. Likewise,
with respect to the two horizontal Lie algebroids, $\mathcal{A}\to A^V$
and $A^H\to M$, the double vector bundle $\mathcal{A}$ is an \LAvb.}

\bigskip

Denote the four anchors by $\tila_V\colon \mathcal{A} \to T(A^H),\
\tila_H\colon \mathcal{A}\to T(A^V),$ and $a_V\colon A^V\to TM,\
a_H\colon A^H\to TM$. As usual we denote all four brackets by $[\ ,\ ]$;
the notation for elements will make clear which structure we are using.

The anchors thus give morphisms of double vector bundles
$$
(\tila_V;\id,a_V;\id)\colon (\mathcal{A};A^H,A^V;M)\to (T(A^H);A^H,TM;M),
$$
$$
(\tila_H;a_H,\id;\id)\colon (\mathcal{A};A^H,A^V;M)\to (T(A^V);TM,A^V;M)
$$
and so define morphisms of their cores; denote these by
$\da_H\colon K\to A^H$ and $\da_V\colon K\to A^V$.

Now return to $\pi^{\#V}$ and Figure~\ref{fig:cofd}.
Since the corner map $\mathcal{A}\to T(A^H)$ is $\tila_V$, the corner map
$A^V\to TM$ is $a_V$ (or $-a_V$ if $Z_H$ is incorporated) and the core
map for the front faces is $\da_H$. Likewise, since the core map of the
left faces is $-\tila^*_V$, the core map of the right faces must be
$-\da_H^*\colon (A^H)^*\to K^*$ (whether or not $Z_H$ is incorporated).
Lastly, the core map of the top faces is
$\pi^\#_V\colon T^*((A^V)^*)\to T((A^V)^*)$, the anchor for the Poisson
structure on $(A^V)^*$ dual to the given Lie algebroid structure on $A^V$.
(These observations are all special cases of \cite[\S3]{Mackenzie:1999}.)

Each of the maps of the core double vector bundles induces on
$T^*M\to (A^V)^*$ the map $-a_V^*$.

Similarly we can analyze $\pi^{\#H}\colon T^*(\mathcal{A}^{*H})\to
T(\mathcal{A}^{*H})$ as a morphism of triple vector bundles, and obtain
$\chigh_H\colon (\mathcal{A}^{*H})^\dagger\to T(K^*)$.

For Condition II, note first that it is automatic that $\tila_V$ is a
morphism of Lie algebroids over $A^H$ and that $a_V$ is a morphism
of Lie algebroids over $M$.

\subsection{Condition II}

{\em The anchors $\tila_V$ and $a_V$ form a morphism of Lie algebroids
$(\tila_V, a_V)$ with respect to the horizontal structure on $\mathcal{A}$
and the prolongation to $TA^H\to TM$ of the structure on $A^H\to M$.
Likewise, the anchors $\tila_H$ and $a_H$ form a morphism of Lie algebroids
with respect to the vertical structure on $\mathcal{A}$ and the prolongation
to $TA^V\to TM$ of the structure on $A^V\to M$.
}

\bigskip

By Condition I and \S\ref{sect:pc},
the Poisson structure on $\mathcal{A}^{*H}\to K^*$ is
linear, and therefore induces a Lie algebroid structure on its dual
$(\mathcal{A}^{*H})^\dagger$. We use $Z_V$ to transfer this to
$\mathcal{A}^{*V}\to K^*$. Similarly the linear Poisson structure on
$\mathcal{A}^{*V}\to K^*$ induces a Lie algebroid structure on
$(\mathcal{A}^{*V})^\dagger\to K^*$.

\subsection{Condition III}

{\em With respect to these structures,
$(\mathcal{A}^{*V}, (\mathcal{A}^{*V})^\dagger)$ is a Lie bialgebroid.
Further, $(\mathcal{A}^{*V}; A^H, K^*; M)$ is an \LAvb\ with respect to the
horizontal Lie algebroid structures and likewise
$(\mathcal{A}^{*H}; K^*, A^V ;M)$
is an \LAvb\ with respect to the vertical structures.
}

\bigskip

\begin{df}                                    \label{df:doubla}
A {\em double Lie algebroid} is a double vector bundle as in
Figure~{\rm \ref{fig:intro}(a)} equipped with Lie algebroid structures on all
four sides such that the above conditions {\rm I, II, III} are satisfied.
\end{df}

For $(\mathcal{A}; A^H, A^V; M)$ a double Lie algebroid, we call
$(\mathcal{A}^{*V}, (\mathcal{A}^{*V})^\dagger)$ the {\em associated Lie
bialgebroid}.

The notion of Lie bialgebroid was defined in \cite{MackenzieX:1994} in
terms of the coboundary operators associated to $A$ and to $A^*$; a more
efficient and elegant reformulation was then given in
\cite{Kosmann-Schwarzbach:1995}. The definition most useful to us here
is quoted below in \ref{theorem:6.2}. For the moment we only need the
following.

Suppose that $(E, E^*)$ is a Lie bialgebroid on base $P$ and denote the
anchors by $e$ and $e_*$. Then we take the Poisson structure on $P$ to
be $\pi^\#_P = e_*\circ e^*$; this is the opposite to \cite{MackenzieX:1994},
but the same as \cite{Kosmann-Schwarzbach:1995}. It follows that $e$ is a
Poisson map (to the tangent lift structure on $TP$) and $e_*$ is
anti--Poisson.

One expects the core of a double Lie algebroid to have a Lie algebroid
structure induced from those on $\mathcal{A}$. However, as \ref{lem:K0}
shows, the straightforward embedding of $K$ in terms of core sections
yields only the zero structure (see also \cite{Mackenzie:DoublaII}). Here
we obtain the correct structure in terms of its dual.

The anchor of $(\mathcal{A}^{*V})^\dagger$ is $\chigh_V$, the appropriate
corner map of the Poisson anchor for $\mathcal{A}^{*V}$. On the other hand,
the anchor for $\mathcal{A}^{*V}$ itself is $\chigh_H\circ Z_V^{-1}.$ So
the Poisson anchor for $K^*$ is
$$
\pi^\#_{K^*} = \chigh_V\circ (Z_V^{-1})^\dagger \circ \chigh_H^\dagger.
$$
Now $Z_V^\dagger = Z_H$ has side map $-\id\colon A^V\to A^V$ and
$\chigh_V$ has side map $a_V$. The side map of $\chigh_H^\dagger$ is the
dual of the core map $-\da_V^*$ of $\chigh_H$. Thus the side map of
$\pi^\#_{K^*}$ is $a_V\circ\da_V$. One likewise checks that the core map
is $-\da_H^*\circ a_H^*$. This proves the first half of the following
result.

\begin{prop}
The anchor $a_K$ for the Lie algebroid structure on $K$ induced by the
Poisson structure on $K^*$ which arises from the Lie bialgebroid structure
on $(\mathcal{A}^{*V}, (\mathcal{A}^{*V})^\dagger)$ is
$a_V\circ\da_V = a_H\circ\da_H$. The maps $\da_H\colon K\to A^H$ and
$\da_V\colon K\to A^V$ are Lie algebroid morphisms.
\end{prop}

\pf
Since $\chigh_V$ is the anchor for $(\mathcal{A}^{*V})^\dagger$, it is
anti--Poisson into $T(K^*)$. Regarding $\chigh_V$ as a morphism of
the right faces in Figure~\ref{fig:cofd}, its core is $-\da_H^*$,
which is therefore anti--Poisson. So $\da^H$ is a morphism of Lie
algebroids.
\boom

The most fundamental example motivating \ref{df:doubla} is for us
the double Lie algebroid of a double Lie groupoid, as constructed in
\cite{Mackenzie:1992}, \cite{Mackenzie:DoublaII}. Some of what is
required in order to verify that the double Lie algebroid of a double
Lie groupoid does satisfy \ref{df:doubla} has been given in
\cite{Mackenzie:1999}, and we recall the details briefly.

In order to proceed, we need to describe the notion of double Lie groupoid
in more detail (see \cite{Mackenzie:1992} and references given there).
A double Lie groupoid consists of a manifold $S$ equipped with two Lie
groupoid structures on bases $H$ and $V$, each of which is a Lie
groupoid over a common base $M$, such that the structure maps (source,
target, multiplication, identity, inversion) of each groupoid structure
on $S$ are morphisms with respect to the other; see Figure~\ref{fig:S}(a).
\begin{figure}[htb]
$$
\begin{matrix}
        &&&&\cr
          &S&\sgpd &V&\cr
          &&&&\cr
          &\vgpd&&\vgpd&\cr
          &&&&\cr
          &H&\sgpd &M&\cr
          &&&&\cr
          &&\mbox{(a)}&&\cr
\end{matrix}
\qquad\qquad
\begin{matrix}
                      &&   &&\cr
                            &A_VS    &\sgpd &AV&\cr
                            &&&&\cr
                            &\Bigg\downarrow&&\Bigg\downarrow&\cr
                            &&&&\cr
                            &H &\sgpd &M&\cr
                            &&&&\cr
                            &&\mbox{(b)}\cr
\end{matrix}
\qquad\qquad
\begin{matrix}
      &&&&\cr
        &A^2S&\mlra   &AV&\cr
        &&&&\cr
        &\Bigg\downarrow& &\Bigg\downarrow&   \cr
        &&&&\cr
        &AH&\mlra &M&\cr
        &&&&\cr
        &&\mbox{(c)}\cr
\end{matrix}
$$
\caption{\ \label{fig:S}}
\end{figure}
One should think of elements of $S$ as squares, the horizontal edges of
which come from $H$, the vertical edges from $V$, and the corner points
from $M$.

Consider a double Lie groupoid $(S;H,V;M)$ as in Figure~\ref{fig:S}(a).
Applying the Lie functor to the vertical structure $S\gpd H$ gives a Lie
algebroid $A_VS\to H$ which has also a groupoid structure over $AV$ obtained
by applying the Lie functor to the structure maps of $S\gpd V$; this is
the {\em vertical \LAgpd}\ of $S$ \cite[\S4]{Mackenzie:1992}, as in
Figure~\ref{fig:S}(b). The Lie algebroid of $A_VS\gpd AV$ is denoted $A^2S$;
there is a double vector bundle structure $(A^2S;AH,AV;M)$ obtained by
applying $A$ to the vector bundle structure of $A_VS\to H$
\cite{Mackenzie:DoublaII}; see Figure~\ref{fig:S}(c). Reversing the order of
these operations, one defines first the {\em horizontal \LAgpd}\
$(A_HS;AH,V;M)$ and then takes the Lie algebroid $A_2S = A(A_HS)$. The
canonical involution $J_S\colon T^2S\to T^2S$ then restricts to an
isomorphism of double vector bundles $j_S\colon A^2S\to A_2S$ and
allows the Lie algebroid structure on $A^2S\to AV$ to be transported to
$A_2S\to AV$. Thus $A_2S$ is a double vector bundle equipped with four
Lie algebroid structures; in \cite{Mackenzie:DoublaII} we called this
the {\em double Lie algebroid of $S$}. The core of both double vector bundles
$A_2S$ and $A^2S$ is $AC$, the Lie algebroid of the core groupoid $C\gpd M$
of $S$ \cite[1.6]{Mackenzie:DoublaII}.

Consider $A_2S$. The structure maps for the horizontal vector bundle
$A_2S\to AV$ are obtained by applying the Lie functor to the structure
maps of $A_HS\to V$ and are therefore Lie algebroid morphisms with
respect to the vertical Lie algebroid structures. The corresponding
statement is true for the vertical vector bundle $A^2S\to AH$ and this
is transported by $j_S$ to $A_2S\to AH$. Thus Condition~I holds.

Let $\ts{a}_V\colon A_2S\to TAH$ denote the anchor for the Lie algebroid of
$A_HS\gpd AH$. Then, as with any Lie groupoid, $\ts{a}_V = A(\ts{\chigh}_V)$
where $\ts{\chigh}_V\colon A_HS\to AH\times AH$ combines the target and
source of $A_HS\gpd AH$. It is easily checked that $\ts{\chigh}_V$ is a
morphism of \LAgpd s over $\chigh_V\colon V\to M\times M$ and $\id\colon
AH\to AH$, and so it follows, by using the methods of
\cite[\S1]{Mackenzie:DoublaII}, that $\ts{a}_V$ is a morphism of Lie algebroids
over $a_V$. Similarly one transports the result for the anchors
$A^2S\to TAV$ and $AH\to TM$. Thus Condition~II is satisfied.

Now consider the bialgebroid condition. The vertical dual
$\mathcal{A}^{*V}$ is $A^*(A_HS)$ and in order to take the dual of this
over $A^*C$ we use the isomorphism ${j'}^V\colon A^*(A_HS)\to A(A^*_VS)$ of
\cite[(21)]{Mackenzie:1999}. This induces $(\mathcal{A}^{*V})^\dagger
\isom A^*(A^*_VS)$.

Now the structure on $\mathcal{A}^{*V}$ itself comes from
$(\mathcal{A}^{*H})^\dagger$. We have $\mathcal{A}^{*H} = A^\sol(A_HS)$ and
the isomorphism $I_H\colon A(A^*_HS)\to A^\sol(A_HS)$ associated to the
\LAgpd\ $A_HS$ in \cite[\S3]{Mackenzie:1999} (see also (\ref{eq:needI})
below) allows us to identify $(\mathcal{A}^{*H})^\dagger$ with
$A^*(A^*_HS)$.

This proves that $(\mathcal{A}^{*V}, (\mathcal{A}^{*V})^\dagger)$ is
isomorphic to $(A^*(A^*_HS), A^*(A^*_VS))$. Now the isomorphism
$\mathcal{D}_H\colon A^*(A^*_HS)\to A(A^*_VS)$ of
\cite[3.9]{Mackenzie:1999} allows this to be written as
$(A(A^*_VS), A^*(A^*_VS))$ and this is the Lie bialgebroid of
$A^*_VS\gpd A^*C$, which was proved to be a Poisson groupoid in
\cite[2.12]{Mackenzie:1999}. Notice that we started with $A_2S$,
defined in terms of the horizontal \LAgpd, and ended with the Lie
bialgebroid of the dual of the vertical \LAgpd.

To make this sketch into a proof, one must ensure that the various
isomorphisms preserve the Poisson structures involved. Rather than do
this, we prove a more general result.

Consider an \LAgpd\ as in Figure~\ref{fig:LAgpd}(a); that is, $\Om$ is
both a Lie algebroid over $G$ and a Lie groupoid over $A$, and each of the
groupoid structure maps is a Lie algebroid morphism; further, the map
$\Om\to A\pback{M} G$ defined by the source and the bundle projection, is a
surjective submersion.
\begin{figure}[htb]
$$
\begin{matrix}
                      &&\tilq   &&\cr
                            &\Om     &\mlra &G &\cr
                            &&&&\cr
                            &\vgpd&&\vgpd&\cr
                            &&&&\cr
                            &A  &\mlra &M&\cr
                            &&q_A&&\cr
                            &&&&\cr
                            &&\mbox{(a)}&&\cr
\end{matrix}
\qquad\qquad\qquad
\begin{matrix}
                      &&A(\tilq) &&\cr
                            &A\Om     &\mlra &AG&\cr
                            &&&&\cr
                   \ts{q}   &\Bigg\downarrow&&\Bigg\downarrow&q_G\cr
                            &&&&\cr
                            &A  &\mlra &M&\cr
                            &&q_A&&\cr
                            &&&&\cr
                            &&\mbox{(b)}&&\cr
\end{matrix}
$$
\caption{\ \label{fig:LAgpd}}
\end{figure}
Applying the Lie functor vertically gives a double vector bundle
$\mathcal{A} = A\Om$ as in Figure~\ref{fig:LAgpd}(b), with Lie algebroid
structures on the vertical sides. It is shown in \cite[\S1]{Mackenzie:DoublaII}
that the Lie algebroid structure of $\Om\to G$ may be prolonged to
$A\Om\to AG$.

That the anchor $\ts{a}\colon A\Om\to TA$ for the Lie algebroid of
$\Om\gpd A$ is a morphism of Lie algebroids over $a_G\colon AG\to TM$
follows as in the case of $A_2S$ above. The anchor ${\bf a}\colon A\Om\to TAG$
for the prolongation structure is $j_G^{-1}\circ A(\tila)$ where
$j_G\colon TAG\to ATG$ is the canonical isomorphism of
\cite[7.1]{MackenzieX:1994}. Since $\tila\colon\Om\to TG$ is a groupoid
morphism over $a\colon A\to TM$, and $j_G$ is an isomorphism of Lie
algebroids over $TM$, it follows that Condition~II is satisfied.
Condition~I is dealt with in the same way.

It was shown in \cite[\S3]{Mackenzie:1999} that $\Om^*\gpd K^*$, the
dual groupoid of $\Om$, together with the Poisson structures on $\Om^*$ and
$K^*$ dual to the Lie algebroid structures on $\Om$ and the core $K$, is a
Poisson groupoid. Thus Condition~III will follow from the next result.

\begin{thm}                               \label{theorem:AOm}
The canonical isomorphism of double vector bundles $R = R^{gpd}$ from
$A^*\Om^*$ to $A^*\Om = \mathcal{A}^{*V}$ is an isomorphism of Lie
bialgebroids
$$
(A^*\Om^*, \Bar{A\Om^*})\to (\mathcal{A}^{*V}, (\mathcal{A}^{*V})^\dagger)
$$
where $(A\Om^*, A^*\Om^*)$ is the Lie bialgebroid of the Poisson
groupoid $\Om^*$ on $K^*$.
\end{thm}

We first recall the map $R$ from \cite[3.8]{Mackenzie:1999}. Associated
with $A\Om$ oriented as in Figure~\ref{fig:LAgpd}(b) there is the pairing
of the vertical and horizontal duals (\ref{eq:3duals}), which we write in
mnemonic form:
$$
\langle A^*\Om, A^\sol\Om\rangle_{K^*} =
\langle A^\sol\Om, A\Om\rangle_{AG} -
\langle A^*\Om, A\Om\rangle_{A}
$$
with the subscripts indicating the bases of the pairings. As in
\cite[(19)]{Mackenzie:1999},
we use the canonical isomorphism $I\colon A\Om^*\to A^\sol\Om$ induced
by the pairing $\llangle\ ,\ \rrangle\colon A\Om^*\pback{AG} A\Om\to\R$
to transfer this to a pairing $\ddagger\ ,\ \ddagger$ of $A^*\Om$ and
$A\Om^*$ over $K^*$ for which
\begin{equation}                                   \label{eq:needI}
\ddagger A^*\Om, A\Om^*\ddagger =
\llangle A\Om^*, A\Om\rrangle -
\langle A^*\Om, A\Om\rangle_{A}.
\end{equation}
Now define $R$ by $\ddagger\mathcal{X}, R(\mathcal{F})\ddagger =
\langle\mathcal{X}, \mathcal{F}\rangle$, where $\mathcal{X}\in A\Om^*,
\mathcal{F}\in A^*\Om^*$, and the pairing on the RHS is the standard one
over $K^*$. We finally arrive at
\begin{equation}                                      \label{eq:3dualsAOm}
\langle\mathcal{X}, \mathcal{F}\rangle_{K^*} +
\langle R(\mathcal{F}), \Xi\rangle_A  =
\llangle \mathcal{X}, \Xi\rrangle
\end{equation}
for compatible $\Xi\in A\Om$. Equivalently,
$Z_V\colon (A^\sol\Om)^\dagger\to A^*\Om$ is given by
\begin{equation}                                  \label{eq:ZRI}
Z_V = R\circ I^\dagger.
\end{equation}

The first part of the following result was stated without proof in
\cite[\S3]{Mackenzie:1999}.

\begin{prop}                             \label{prop:RI}
{\em (i)} The map $R\colon A^*\Om^*\to A^*\Om$ is anti--Poisson from
the Poisson structure dual to the Lie algebroid of $\Om^*\gpd K^*$ to
the Poisson structure dual to the Lie algebroid of $\Om\gpd A$.

{\em (ii)} The map $I\colon A\Om^*\to A^\sol\Om$ is Poisson from the Poisson
structure induced on $A\Om^*$ {\em \cite{Weinstein:1988}} by the Poisson
groupoid structure on $\Om^*\gpd K^*$, to the Poisson structure dual to
the prolonged Lie algebroid structure on $A\Om\to AG$.
\end{prop}

\pf
It suffices \cite{Weinstein:1988} to prove that the graph of $R$ is
coisotropic in $A^*\Om^*\times A^*\Om$. Let $S = \Om^*\pback{G}\Om$
and write $F\colon S\to\R$ for the pairing. Then $F$ is a groupoid
morphism, where $S$ is the pullback groupoid over $K^*\pback{M}A$,
and so, as in \cite[3.7]{Mackenzie:1999}, we can apply the Lie
functor and get $\llangle\ ,\ \rrangle = A(F)\colon AS\to\R$. This is
linear and so defines a section $\nu$ of the dual of $AS$, which is closed
since $A(F)$ is a morphism. So by \cite[4.6]{MackenzieX:1998}, the image
of $\nu$ is coisotropic.

It remains to show that the image of $\nu$ coincides with the graph
of $R$. The image of $\nu$ consists of those $(\mathcal{F}, \mathcal{X})\in
A^*_\kappa\Om\times A^*_Y\Om$ such that
$$
\langle(\mathcal{F}, \Phi\rangle, (\mathcal{X}, \Xi)\rangle =
A(F)(\mathcal{X}, Xi)
$$
for all $(\mathcal{X}, \Xi)\in AS$ compatible with $(\mathcal{F}, \Phi)$.
As in \cite[5.5]{MackenzieX:1994}, this equation expands to
(\ref{eq:3dualsAOm}).

We leave the proof of (ii) to the reader.
\boom

\noindent
{\sc Proof of Theorem \ref{theorem:AOm}:} We must first show that $R$ is an
isomorphism of Lie algebroids $A^*\Om^*\to \mathcal{A}^{*V}$. Now the Lie
algebroid structure on $\mathcal{A}^{*V}$ is induced from
$(\mathcal{A}^{*H})^\dagger$ via $Z_V$. So what we have to show is that
$Z_V^{-1}\circ R\colon A^*\Om^*\to (A^\sol\Om)^\dagger$ is an isomorphism
of Lie algebroids, and this is equivalent to the dual over $K^*$ being
Poisson. This dual is, using (\ref{eq:ZRI}), $I^{-1}\colon A^\sol\Om\to
A\Om^*$, and so the result follows from \ref{prop:RI}(ii) above.

Secondly we must show that
$R^\dagger\colon (\mathcal{A}^{*V})^\dagger\to\Bar{A\Om^*}$ is an isomorphism
of Lie algebroids over $K^*$. (Note that the minus sign is in the
bundle over $K^*$.) This is equivalent to showing that the dual
$R\colon\Bar{A^*\Om^*}\to A^*\Om$ is Poisson, and this is
\ref{prop:RI}(i) above.
\boom

In summary, we have proved:

\begin{thm}                                    \label{thm:infinv}
The double Lie algebroid $(A\Om;A,AG;M)$ of an \LAgpd\
$(\Om;A,G;M)$, as constructed in {\em \cite[\S2]{Mackenzie:DoublaII}}, is a
double Lie algebroid as defined in {\em \ref{df:doubla}}.

In particular, the double Lie algebroids $(A_2S;AH,AV;M)$ and
$(A^2S;AH,AV;M)$ of a double Lie groupoid $(S;H,V;M)$, as constructed
in {\em \cite[\S3]{Mackenzie:DoublaII}}, are double Lie algebroids as
defined in {\em \ref{df:doubla}}.
\end{thm}

\begin{ex}\rm                                          \label{ex:TA}
Let $A$ be any Lie algebroid on $M$. Then
$\Om = A\times A$ has an \LAgpd\ structure over $M\times M$ and $A$,
and the associated double Lie algebroid constructed in
\cite[\S1]{Mackenzie:DoublaII} is $(TA;A,TM;M)$.

The associated duals are $\mathcal{A}^{*V} = T^*A$ and
$\mathcal{A}^{*H} = T^\sol A$. Using $R$ and $I$ as in \cite{MackenzieX:1994},
these can be identified with $T^*A^*$ and $T(A^*)$, as bundles over $A^*$.
The Lie algebroid structure on $T^*A^*$ is the cotangent of the dual Poisson
structure on $A^*$. The Lie algebroid structure on $T(A^*)$ is the standard
tangent bundle structure. The Lie bialgebroid associated to $TA$ is the
standard Lie bialgebroid $(T^*P, TP)$ for $P = A^*$.
\end{ex}

\begin{ex}\rm
Taking $A = TM$ in the previous example, we see that $T^2M$ is a double
Lie algebroid, the associated bialgebroid of which is $(T^*T^*M,TT^*M)$.
This is a Lie bialgebroid over $T^*M$, the induced Poisson structure
being the standard symplectic structure.
\end{ex}

The double Lie algebroids considered in the next two sections do not
necessarily have an underlying \LAgpd.

\section{THE DOUBLE LIE ALGEBROID OF A LIE \break BIALGEBROID}
\label{sect:dlalba}

Here we use the following criterion for a Lie bialgebroid.

\begin{thm}                                        \label{theorem:6.2}
{\bf \cite[6.2]{MackenzieX:1994}}
Let $A$ be a Lie algebroid on $M$ such that its dual vector bundle $A^{*}$
also has a Lie algebroid structure. Denote their anchors by $a, a_*$.
Then $(A, A^{*})$ is a Lie bialgebroid if and only if
\begin{equation}                         \label{main}
T^*(A^*)\buildrel{R}\over\longrightarrow T^*(A)
\buildrel\pi^\#_A\over\longrightarrow TA
\end{equation}
is a Lie algebroid morphism over $a_*$, where the domain
$T^{*}(A^{*})\to A^{*}$ is the cotangent Lie algebroid induced by the
Poisson structure on $A^{*}$, and the target $TA\to TM$ is the tangent
prolongation of $A$.
\end{thm}

Consider a Lie algebroid $A$ on $M$ together with a Lie algebroid
structure on the dual, not {\em a priori} related to that on $A$. The
structure on $A^*$ induces a Poisson structure on $A$, and this gives
rise to a cotangent Lie algebroid $T^*A\to A$. Equally, the Lie algebroid
structure on $A$ induces a Poisson structure on $A^*$ and this gives rise
to a cotangent Lie algebroid $T^*A^*\to A^*$. We transfer this latter
structure to $T^*A\to A^*$ via $R$.

There are now four Lie algebroid structures on the four sides of
$\mathcal{A} = T^*A$ as in Figure~\ref{fig:intro}(c).

\begin{thm}                                       \label{theorem:Manin}
Let $A$ be a Lie algebroid on $M$ such that its dual vector bundle $A^{*}$
also has a Lie algebroid structure. Then $(A, A^*)$ is a Lie bialgebroid
if and only if $\mathcal{A} = T^*A$, with the structures just described, is a
double Lie algebroid.
\end{thm}

\pf
Assume that $(A, A^*)$ is a Lie bialgebroid. The vertical structure
on $\mathcal{A}$ is the cotangent Lie algebroid structure for the Poisson
structure on $A$. The anchor of this is a morphism of double vector
bundles $\pi^\#_A\colon T^*A\to TA$ over $a_*\colon A^*\to TM$ and $\id_A$,
inducing $-a^*_*\colon T^*M\to A$ on the cores. Now the horizontal structure
has the cotangent Lie algebroid structure for the Poisson structure on
$A^*$, transported via $R = R_A\colon T^*A^*\to T^*A$. So the condition that
$\pi^\#_A$ is a morphism of Lie algebroids over $a_*$ with respect to the
horizontal structure is precisely \ref{theorem:6.2}.

On the other hand, the anchor for the horizontal structure is
$$
\pi^\#_{A^*}\circ R^{-1}\colon T^*A\to T(A^*),
$$
and this is a morphism of double vector bundles over $a\colon A\to TM$
and $\id_{A^*}$, inducing $+a^*$ on the cores. Since $R^{-1} = R_{A^*}$,
the condition that this anchor be a morphism with respect to the
vertical structure is precisely the dual form of \ref{theorem:6.2}, to which
\ref{theorem:6.2} is equivalent by \cite[3.10]{MackenzieX:1994} or
\cite{Kosmann-Schwarzbach:1995}.

The vertical dual of $\mathcal{A}$ is the tangent double vector bundle as in
Figure~\ref{fig:intro}(b). Being the dual of a Lie algebroid, the vertical
structure of this has a Poisson structure; this is the tangent lift
of the Poisson structure on $A$ \cite[5.6]{MackenzieX:1994}. The
corresponding Poisson tensor is
$$
\pi^\#_{TA} = J_A\circ T(\pi^\#_A)\circ\theta_A^{-1}
$$
where $J_A\colon T^2A\to T^2A$ is the canonical involution for the
manifold $A$ and
$\theta_A\colon T(T^*A)\to T^*(TA)$ is the canonical map $\alpha$ of
\cite{Tulczyjew}, denoted $J'$ in \cite[5.4]{MackenzieX:1994}.

We must check that this Poisson structure coincides with that induced from
the dual of $\mathcal{A}^{*H}$. Note that we have $K^* = TM$ here and to avoid
confusion we drop the ${}^\dagger$ notation in this case and denote all
duals over $TM$ by ${}^\sol$. Duals over $A^*$ will be denoted
${}^{*A^*}$.

Consider $I\circ R^{*A^*}\colon (T^*A)^{*A^*}\to T^\sol A$. This preserves
the core $A^*$ and the side $A^*$ but reverses the side $TM$. Define
$W = -I\circ R^{*A^*}$ where the minus is for the bundle over $A^*$.
Then $W^\sol\colon TA\to (\mathcal{A}^{*H})^\sol$ and the reader can
check that $Z_V^{-1} = \dminus W^\sol$ where the heavy minus is over $TM$.
See alternatively \cite[3.3]{Mackenzie:1999}.

To prove that $Z_V\colon (\mathcal{A}^{*V})^\sol \to TA$ is an isomorphism
of Lie algebroids over $TM$ we must show that $W$ is an anti--Poisson
map. This may be done directly or by observing that $-W$ is, in terms
of the double Lie algebroid $\mathcal{A}' = TA$ of \ref{ex:TA}, the map
$(Z'_V)^\dagger = Z'_H$.

So we have $\mathcal{A}^{*V} = TA$ and $(\mathcal{A}^{*V})^\sol = T^\sol A$ and
Condition~III follows from the next result. We use $I$ to replace $T^\sol A$
by $TA^*$.

\begin{lem}
Given that $(A,A^*)$ is a Lie bialgebroid on $M$, the tangent prolongation
structures make $(TA, TA^*)$ a Lie bialgebroid on $TM$ with respect to the
tangent pairing.
\end{lem}

\pf
We use the bialgebroid criterion of \ref{theorem:6.2}. We must
prove that
\begin{equation}                                \label{eq:comp}
T^*(TA^*)\buildrel\Tilde{R}\over\longrightarrow T^*(TA)
\buildrel\pi^\#_{TA}\over\longrightarrow T^2A
\end{equation}
is a morphism of Lie algebroids over the anchor of $TA^*$; by
\cite[5.1]{MackenzieX:1994} this anchor is
$J_M\circ T(a_*)\colon TA^*\to T^2M$.
Here $\Tilde{R}$ is the canonical map $R$ for $TA\to TM$, transported
using $I\colon TA^*\to T^\sol A$. The domain of (\ref{eq:comp}) is the
cotangent Lie algebroid for the Poisson structure on $TA^*$, which Poisson
structure---again by \cite[5.6]{MackenzieX:1994}---is both the tangent lift
of the Poisson structure on $A^*$ and the dual (via $I$) of the prolongation
Lie algebroid structure on $TA$. The target of (\ref{eq:comp}) is the
iterated tangent prolongation of the Lie algebroid structure of $A$.

Now $\Tilde{R} = \theta_A\circ T(R_A)\circ\theta_{A^*}^{-1}$ and so
$$
\pi^\#_{TA}\circ \Tilde{R} =
J_A\circ T(\pi^\#_A\circ R_A)\circ\theta_{A^*}^{-1}.
$$
We know that $\pi^\#_A\circ R_A\colon T^*A^*\to TA$ is a morphism of
Lie algebroids over $a_*$, so $T(\pi^\#_A\circ R_A)$ is a morphism of
the prolongation structures over $T(a_*)$. We need two further
observations.

Firstly, for any Poisson manifold, $\theta_P\colon T(T^*P)\to T^*(TP)$
is an isomorphism of Lie algebroids over $TP$ from the tangent prolongation
of the cotangent Lie algebroid structure on $T^*P$ to the cotangent Lie
algebroid of the tangent Poisson structure \cite[2.13]{Mackenzie:DoublaII}.
We apply this to $P = A^*$.

Secondly, $J_A\colon T^2A\to T^2A$ is a Lie algebroid automorphism over
$J_M$ of the iterated prolongation of the given Lie algebroid
structure on $A$.

Putting these facts together, we have that $\pi^\#_{TA}\circ\Tilde{R}$ is a
Lie algebroid morphism.
\boom

Now conversely suppose that $A$ is a Lie algebroid on $M$ and that $A^*$
has a Lie algebroid structure, not {\em a priori} related to the structure
on $A$. Consider $\mathcal{A} = T^*A$ with the two cotangent Lie algebroid
structures arising from the Poisson structures on $A^*$ and $A$, and
suppose that these structures make $\mathcal{A}$ a double Lie algebroid.

Then in particular the anchor of the horizontal structure must be a Lie
algebroid morphism with respect to the other structures, as in Condition~II,
and this is
$$
T^*A\buildrel{R_{A^*}}\over\longrightarrow T^*A^*
\buildrel\pi^\#_{A^*}\over\longrightarrow TA^*
$$
That this be a Lie algebroid morphism over $a\colon A\to TM$ is
precisely the dual form of \ref{theorem:6.2}.

This completes the proof of Theorem \ref{theorem:Manin}.
\boom

Recall the Manin triple characterization of a Lie bialgebra, as given
in \cite{LuW:1990}: Given a Lie bialgebra $({\mathfrak g}, {\mathfrak g}^*)$
the vector space direct sum ${\mathfrak d} = {\mathfrak g}\oplus{\mathfrak g}^*$
has a Lie algebra bracket defined in terms of the two coadjoint
representations. This bracket is invariant under the pairing
$\langle X + \phi, Y + \psi\rangle = \psi(X) + \phi(Y)$ and both ${\mathfrak g}$
and ${\mathfrak g}^*$ are coisotropic subalgebras. Conversely, if a Lie
algebra ${\mathfrak d}$ is a vector space direct sum ${\mathfrak g}\oplus{\mathfrak h}$,
both of which are coisotropic with respect to an invariant pairing
of ${\mathfrak d}$ with itself, then ${\mathfrak h}\isom{\mathfrak g}^*$ and
$({\mathfrak g}, {\mathfrak h})$ is a Lie bialgebra, with ${\mathfrak d}$ as the
double.

Two aspects of this result concern us here. Firstly, it characterizes the
notion of Lie bialgebra in terms of a single Lie algebra structure on
${\mathfrak d}$, the conditions being expressed in terms of the simple notion
of pairing. Secondly, the roles of the two Lie algebras ${\mathfrak g}$ and
${\mathfrak g}^*$ are completely symmetric; it is an immediate consequence of
the Manin triple result that $({\mathfrak g}, {\mathfrak g}^*)$ is a Lie bialgebra
if and only if $({\mathfrak g}^*, {\mathfrak g})$ is so.

In considering a corresponding characterization for Lie bialgebroids, the
most important difference to note is that the structure on the double is
no longer over the same base as the given Lie algebroids. This is to be
expected in view of the results of \cite[\S2]{Mackenzie:1992} for the
double groupoid case. There it is proved that if a double groupoid
$(S;H,V;M)$ has trivial core (that is, the only elements of $S$ to have
two touching sides which are identity elements, are those which are
identities for both structures), then there is a third groupoid structure
on $S$, over base $M$, called in \cite[p.200]{Mackenzie:1992} the
{\em diagonal structure} and denoted $S_D$. With respect to this structure
on $S$, the identity maps from $H$ and $V$ into $S$ are morphisms over
$M$, and $S$ as a manifold is $H\pback{M}V$. This diagonal structure is, in
the case where $H$ and $V$ are dual Poisson groups, precisely the structure
of the double group. The existence of the diagonal structure, in the general
formulation given in \cite[\S2]{Mackenzie:1992}, depends crucially on the
fact that $S$ has trivial core; that is, that $S$ is vacant. Since the core
of the double vector bundle $T^*A$, for $A$ a vector bundle on $M$, is
$T^*M$, we expect $T^*A$ to possess a Lie algebroid structure over $M$ only
when $M$ is a point.

The role played in the bialgebra case by the Lie algebra structure of
${\mathfrak d} = {\mathfrak g}\oplus{\mathfrak g}^*$ is taken, for Lie bialgebroids,
by the two structures on $T^*A$ (the bases of which are $A$ and $A^*$). In place
of a characterization in terms of a single Lie algebroid structure on
$T^*A$ with base $M$, Theorem \ref{theorem:Manin} gives a characterization in
terms of the two Lie algebroid structures on $T^*A$. The role of the pairing
in the bialgebra case is taken in \ref{theorem:Manin} by the isomorphism $R$.

There are considerable differences between the general notion of Courant
algebroid and that of double Lie algebroid. Firstly, it seems clear that
both concepts include examples not covered by the other. Open Question~1
of \cite{LiuWX:1997} asks for Courant algebroids which do not arise from
Lie bialgebroids; it will be very surprising if these do not exist. There
are certainly many classes of double Lie algebroids which are not the
cotangent doubles of Lie bialgebroids, and it is not clear whether the
correspondence between the Courant algebroid of a Lie bialgebroid and its
cotangent double Lie algebroid can be extended to more general cases of
either concept.

The notion of Courant algebroid provides a generalization to arbitrary
bialgebroids of the Manin triple theorem and, perhaps most importantly,
provides a setting for the study of Dirac structures.
The notion of double Lie algebroid, on the other hand, while
providing an alternative form of the Manin triple theorem for Lie
bialgebroids, is also intended to give (as is argued in \S\ref{sect:adla},
\S\ref{sect:mpvdla}, and the Introduction to \cite{Mackenzie:1992}) a
unification of iterated and second--order constructions in the foundations
of differential geometry, and a general setting for duality phenomena.
The two concepts are not only different, but differ in the nature of
the problems which they are designed to address.

\section{MATCHED PAIRS AND VACANT DOUBLE LIE \break ALGEBROIDS}
\label{sect:mpvdla}

The history of matched pairs of Lie algebras was briefly summarized in
the Introduction. The corresponding concept of matched pair of Lie groups
\cite{LuW:1990}, \cite{Majid:1990} was extended to groupoids
in \cite{Mackenzie:1992}. In \cite{Mokri:1997}, Mokri differentiated the
twisted automorphism equations of \cite{Mackenzie:1992} to obtain conditions
on a pair of Lie algebroid representations, of $A$ on $B$ and of $B$ on $A$,
which ensure that the direct sum vector bundle $A\oplus B$ has a Lie
algebroid structure with $A$ and $B$ as subalgebroids. We quote the
following.

\begin{df}  {\bf \cite[4.2]{Mokri:1997}}   \label{df:mokri}
Let $A$ and $B$ be Lie algebroids on base $M$, with anchors $a$ and $b$,
and let $\rho\colon A\to\CDO(B)$ and $\sigma\colon B\to\CDO(A)$ be
representations of $A$ on the vector bundle $B$ and of $B$ on the vector
bundle $A$. Then $A$ and $B$ together with $\rho$ and $\sigma$ form a
{\em matched pair} if the following equations hold for all $X, X_1, X_2
\in \Ga A,\ Y, Y_1, Y_2\in\Ga B$:
\begin{gather*}
\rho_X([Y_1, Y_2]) = [\rho_X(Y_1), Y_2] + [Y_1, \rho_X(Y_2)]
    +\rho_{\sigma_{Y_2}(X)}(Y_1) -  \rho_{\sigma_{Y_1}(X)}(Y_2),\\
\sigma_Y([X_1, X_2]) = [\sigma_Y(X_1), X_2] + [X_1, \sigma_Y(X_2)]
    +\sigma_{\rho_{X_2}(Y)}(X_1) -  \sigma_{\rho_{X_1}(Y)}(X_2),\\
a(\sigma_Y(X)) - b(\rho_X(Y)) = [b(Y), a(X)].
\end{gather*}
\end{df}

The notation $\CDO(E)$ is defined in \S\ref{sect:pc}.

\begin{prop} {\bf \cite[4.3]{Mokri:1997}}        \label{prop:mokri}
Given a matched pair, there is a Lie algebroid structure on the direct sum
vector bundle $A\oplus B$, with anchor $c(X\oplus Y) = a(X) + b(Y)$ and
bracket
\begin{eqnarray*}
\lefteqn{[X_1\oplus Y_1, X_2\oplus Y_2] =} \\
& & \{[X_1, X_2] + \sigma_{Y_1}(X_2) - \sigma_{Y_2}(X_1)\}\oplus
    \{[Y_1, Y_2] + \rho_{X_1}(Y_2) - \rho_{X_2}(Y_1)\}.
\end{eqnarray*}
Conversely, if $A\oplus B$ has a Lie algebroid structure making
$A\oplus 0$ and $0\oplus B$ Lie subalgebroids, then $\rho$ and $\sigma$
defined by $[X\oplus 0, 0\oplus Y] = -\sigma_Y(X) \oplus \rho_X(Y)$
form a matched pair.
\end{prop}

We now show that matched pairs correspond precisely to double Lie algebroids
with zero core. The following definition is a natural sequel to
\cite[2.11, 4.10]{Mackenzie:1992}.

\begin{df}
A double Lie algebroid $(\mathcal{A}; A^H, A^V; M)$ is {\em vacant} if the
combination of the two projections,
$(\tilq_V, \tilq_H)\colon \mathcal{A}\to A^H\pback{M} A^V$ is a diffeomorphism.
\end{df}

Consider a vacant double Lie algebroid, which we will write here as
$(\mathcal{A}; A, B; M)$. Note that $\mathcal{A}\to A^H$ and
$\mathcal{A}\to A^V$ are the pullback bundles $q_A^!B$ and $q_B^!A$. The two
duals are $\mathcal{A}^{*H} = A^*\oplus B$ and
$\mathcal{A}^{*V} = A\oplus B^*$, as vector bundles over $M$, and the duality
is (see \cite[3.4]{Mackenzie:1999})
\begin{equation}                                         \label{eq:pairing}
\langle X + \psi, \phi + Y\rangle =
                           \langle \phi, X \rangle - \langle\psi, Y\rangle.
\end{equation}

The horizontal bundle projection $\tilq_A\colon \mathcal{A}\to B$ is a
morphism of Lie algebroids over $q_A\colon A\to M$ and since it is a
fibrewise surjection, it defines an action of $B$ on $q_A$ as in
\cite[\S2]{HigginsM:1990a}. Namely, each section $Y$ of $B$ induces the
pullback section $1\otimes Y$ of $q_A^!B$ and this induces a vector field
$\eta(Y) = \tilb(1\otimes Y)$ on $A$, where $\tilb\colon\mathcal{A}\to TA$ is
the anchor of the vertical structure. By Conditions I and II, $\eta(Y)$ is
linear over the vector field $b(Y)$ on $M$, in the sense of
\cite[\S1]{MackenzieX:1998}; that is, $\eta(Y)$ is a vector bundle morphism
$A\to TA$ over $b(Y)\colon M\to TM$. It follows that $\eta(Y)$ defines
covariant differential operators $\sigma^{(*)}_Y$ on $A^*$ and $\sigma_Y$
on $A$ by
\begin{equation}                          \label{eq:sigma}
\eta(Y)(\ell_\phi) = \ell_{\sigma^{(*)}_Y(\phi)},\quad
\langle\phi, \sigma_Y(X)\rangle = b(Y)\langle\phi, X\rangle
- \langle\sigma^{(*)}_Y(\phi), X\rangle
\end{equation}
where $\phi\in\Ga A^*,\ X\in\Ga A,$ and $\ell_\phi$ denotes the function
$A\to\R,\ X\mapsto\langle\phi(q_AX), X\rangle$; see
\cite[\S2]{MackenzieX:1998}. Since $\tilq_A$ is a Lie algebroid
morphism, it follows that $\sigma$ is a representation of $B$ on the
vector bundle $A$.

Likewise, $\tilq_B$ is a morphism of Lie algebroids over $q_B$ and for
each $X\in\Ga A$ we obtain a linear vector field $\xi(X)\in \mathcal{X}(B)$
over $a(X)$. We define covariant differential operators
$\rho^{(*)}_X$ on $B^*$ and $\rho_X$ on $B$ by
\begin{equation}                             \label{eq:rho}
\xi(X)(\ell_\psi) = \ell_{\rho^{(*)}_X(\psi)},\quad
\langle\psi, \rho_X(Y)\rangle = a(X)\langle\psi, Y\rangle
- \langle\rho^{(*)}_X(\psi), Y\rangle.
\end{equation}
Again, $\rho^{(*)}_X$ and $\rho_X$ are representations of $A$.

In fact (see \cite[\S2]{HigginsM:1990a}) the two Lie algebroid structures
on $\mathcal{A}$ are action Lie algebroids determined by the actions
$Y\mapsto\eta(Y)$ and $X\mapsto\xi(X)$. It follows that the dual Poisson
structures are semi--direct in a general sense, but we will proceed
on an ad hoc basis.

For a general vector bundle, the functions on the dual are generated by
the linear functions and the pullbacks from the base manifold. In the case
of a pullback bundle such as $q_B^!A$, one can refine this description a
little further.

Since the two vector bundle structures on $A\times_M B$ are pullbacks the
four classes of functions used in \S\ref{sect:pc} simplify slightly.
Namely, if $\pi_B\colon q_B^!A^*\to A^*$ is
$(\phi, Y)\mapsto \phi$, and $\tilq_{A^*}\colon q_B^!A^*\to B$ is the bundle
projection, then the functions on $q_B^!A^*$ are generated by all
$$
\ell_X\circ\pi_B,\qquad
\ell_\psi\circ\tilq_{A^*}\qquad\mbox{and}\qquad
f\circ q_B\circ\tilq_{A^*}
$$
where $X\in\Ga A,\ \psi\in\Ga B^*$ and $f\in C(M)$. Now the Poisson
structure on $q_B^!A^*$ is characterized by
\begin{gather}
\{\ell_{X_1}\circ\pi_B,\ \ell_{X_2}\circ\pi_B\}
= \ell_{[X_1,X_2]}\circ\pi_B,\nonumber \\
\{\ell_{X}\circ\pi_B,\ \ell_\psi\circ\tilq_{A^*}\}
= \ell_{\rho^{(*)}_X(\psi)}\circ \tilq_{A^*}, \label{eq:poissA}  \\
\{\ell_{X}\circ\pi_B,\ f\circ q_B\circ\tilq_{A^*}\}
= a(X)(f)\circ q_B\circ \tilq_{A^*},\nonumber \\
\{F_1\circ\tilq_{A^*},\ F_2\circ\tilq_{A^*}\}
= 0,\nonumber
\end{gather}
where $F_1, F_2$ are any smooth functions on $B$. Similarly,
\begin{gather}
\{\ell_{Y_1}\circ\pi_A,\ \ell_{Y_2}\circ\pi_A\}
= \ell_{[Y_1,Y_2]}\circ\pi_A,\nonumber \\
\{\ell_{Y}\circ\pi_A,\ \ell_\phi\circ\tilq_{B^*}\}
= \ell_{\sigma^{(*)}_Y(\phi)}\circ \tilq_{B^*},\nonumber \\
\{\ell_{Y}\circ\pi_A,\ f\circ q_A\circ\tilq_{B^*}\}
= b(Y)(f)\circ q_A\circ \tilq_{B^*}, \label{eq:poissB}  \\
\{G_1\circ\tilq_{B^*},\ G_2\circ\tilq_{B^*}\}
= 0,\nonumber
\end{gather}

Now these Poisson structures induce Lie algebroid structures on the
direct sum bundles $A\oplus B^*$ and $A^*\oplus B$ over $M$. Consider
first a section $X\oplus\psi$ of $\Ga(A\oplus B^*)$. Via the pairing
(\ref{eq:pairing}), this induces a linear function on $q_B^!A^*$, namely
$$
\ell^\dagger_{X\oplus\psi} = \ell_X\circ\pi_B - \ell_\psi\circ\tilq_{A^*}
$$
where $\ell^\dagger$ refers to the pairing (\ref{eq:pairing}). By following
through the equations (\ref{eq:poissA}) and (\ref{eq:poissB}) one
obtains the following.

\begin{lem}                                  \label{lem:sdps}
The Lie algebroid structure on $A\oplus B^*$ induced as above has anchor
$e(X\oplus\psi) = a(X)$ and bracket
\begin{equation}                            \label{eq:sdprho}
[X_1\oplus\psi_1, X_2\oplus\psi_2] =
   [X_1, X_2]\oplus\{\rho^{(*)}_{X_1}(\psi_2) - \rho^{(*)}_{X_2}(\psi_1)\}.
\end{equation}

The Lie algebroid structure on $A^*\oplus B$ induced as above has anchor
$e_*(\phi\oplus Y) = -b(Y)$ and bracket
\begin{equation}                              \label{eq:sdpsigma}
[\phi_1\oplus Y_1, \phi_2\oplus Y_2] =
   \{\sigma^{(*)}_{Y_2}(\phi_1) - \sigma^{(*)}_{Y_1}(\phi_2)\}\oplus[Y_2, Y_1].
\end{equation}
\end{lem}

Thus $A\oplus B^*$ is the semi--direct product (over the base $M$, in the
sense of \cite{Mackenzie:LGLADG}) of $A$ with the vector bundle $B^*$
with respect to $\rho^{(*)}$. However $A^*\oplus B$ is the opposite Lie
algebroid to the semi--direct product of $B$ with $A^*$.

We can now apply Condition III to $A\oplus B^*$ and $A^*\oplus B$. For
brevity write $E = A\oplus B^*$. Recall \cite{MackenzieX:1994},
\cite{Kosmann-Schwarzbach:1995} that $(E, E^*)$ is a Lie bialgebroid if
and only if
\begin{equation}                         \label{eq:full}
d^E[\phi_1 \oplus Y_1, \phi_2\oplus Y_2] =
   [d^E(\phi_1\oplus Y_1), \phi_2\oplus Y_2] +
   [\phi_1\oplus Y_1, d^E(\phi_2\oplus Y_2)]
\end{equation}
for all $\phi_1\oplus Y_1, \phi_2\oplus Y_2\in \Ga E^*$. It follows
\cite[3.4]{MackenzieX:1994} that for any $f\in C(M),\ X\oplus\psi\in\Ga E$,
\begin{equation}                     \label{eq:lied}
L_{d^Ef}(X\oplus\psi) + [d^{E^*}f, X\oplus\psi] = 0.
\end{equation}
It is easy to check that
$$
d^Ef = d^Af\oplus 0,\qquad d^{E^*}f = 0\oplus d^Bf;
$$
note that these imply that the Poisson structure induced on $M$ by the Lie
bialgebroid $(E, E^*)$ is zero. Now the Lie derivative in (\ref{eq:lied})
is a standard Lie derivative for $E^*$ and so
\begin{eqnarray}                        \label{eq:shift1}
\langle L_{d^Ef}(X\oplus\psi), \phi\oplus Y\rangle
   & = & e_*(d^Ef)\langle X\oplus\psi, \phi\oplus Y\rangle
         -\langle X\oplus\psi, [d^Ef, \phi\oplus Y]\rangle \nonumber\\
   & = & 0 - \langle X\oplus\psi, \sigma^{(*)}_Y(d^Af)\oplus 0\rangle \\
   & = & \langle d^Af, \sigma_Y(X)\rangle - b(Y)\langle d^Af, X\rangle \nonumber\\
   & = & a(\sigma_Y(X))(f) - b(Y)a(X)(f)\nonumber
\end{eqnarray}
where we used (\ref{eq:sigma}) and (\ref{eq:sdpsigma}). Expanding out the
bracket term in (\ref{eq:lied}) in a similar way, we obtain the third
equation in \ref{df:mokri}.

Now consider the bialgebroid equation (\ref{eq:full}). We calculate this
in the case $\phi_1 = \phi_2 = 0$, with arguments $0\oplus \psi_1,\
X_2\oplus 0.$ With these values we refer to (\ref{eq:full}) as equation
(\ref{eq:full}$)_0$. First we need the following lemma, which is a
straightforward calculation.

\begin{lem}
$$
d^E(\phi\oplus Y)(X_1\oplus\psi_1, X_2\oplus\psi_2) =
   (d^A\phi)(X_1, X_2) + \langle\psi_1, \rho_{X_2}(Y)\rangle
            - \langle\psi_2, \rho_{X_1}(Y)\rangle.
$$
\end{lem}

The LHS of (\ref{eq:full}$)_0$ is easily seen to be
$\langle\psi_1, \rho_{X_2}[Y_2, Y_1]\rangle.$ On the RHS, consider
the second term first. Regarding the bracket as a Lie derivative, we have
\begin{eqnarray}                    \label{eq:shift2}
\langle L_{0\oplus Y_1}(d^E(0\oplus Y_2)),
   (0\oplus\psi_1)\wedge(X_2\oplus 0)\rangle
& = & L_{0\oplus Y_1}\langle d^E(0\oplus Y_2), (0\oplus\psi_1)
        \wedge(X_2\oplus 0)\rangle\\
& & \mbox{} - \langle d^E(0\oplus Y_2), L_{0\oplus Y_1}((0\oplus\psi_1)
        \wedge(X_2\oplus 0))\rangle
\nonumber
\end{eqnarray}
Since $e_*(0\oplus Y_1) = -b(Y_1)$, the first term is
$-b(Y_1)\langle\psi_1, \rho_{X_2}(Y_2)\rangle.$ For the second term we
need the following lemma.

\begin{lem}                              \label{lem:lied}
For any $\phi\oplus Y\in \Ga E^*$ and $X\oplus\psi\in\Ga E$, we have
$L_{\phi\oplus Y}(X\oplus\psi) = -\sigma_Y(X)\oplus \Bar{\psi}$ where
for any $Y'\in B$,
$$
\langle\Bar{\psi}, Y'\rangle = -b(Y)\langle\psi, Y'\rangle
   +\langle\sigma^{(*)}_{Y'}(\phi), X\rangle + \langle\psi, [Y,Y']\rangle.
$$
\end{lem}

\pf
This is a Lie derivative for $E^*$ and applying the same device as in
(\ref{eq:shift1}) we have, for any $\phi\oplus Y'\in\Ga E^*$,
$$
\langle\phi'\oplus Y', L_{\phi\oplus Y}(X\oplus\psi)\rangle
= \mbox{} -\langle\phi', \sigma_Y(X)\rangle + b(Y)\langle\psi, Y'\rangle
-\langle\sigma^{(*)}_{Y'}(\phi), X\rangle + \langle\psi, [Y',Y]\rangle.
$$
Setting $Y' = 0$ and $\phi' = 0$ in turn gives the result.
\boom

Now expand out the Lie derivative of the wedge product in (\ref{eq:shift2})
and apply Lemma \ref{lem:lied}. One obtains for the second term on the
RHS of (\ref{eq:full}$)_0$
$$
-\langle\psi_1, [Y_1, \rho_{X_2}(Y_2)]\rangle
+ \langle\psi_1, \rho_{\sigma_{Y_1}(X_2)}(Y_2)\rangle.
$$
The first term on the RHS of (\ref{eq:full}$)_0$ is easily obtained from
this, and combining with the LHS we have the first equation in
\ref{df:mokri}. The second equation is obtained in a similar way from the
dual form of (\ref{eq:full}).

This completes the proof of the first part of the following result.

\begin{thm}                                     \label{theorem:mp}
Let $(\mathcal{A};A,B;M)$ be a vacant double Lie algebroid. Then the two
Lie algebroid structures on $\mathcal{A}$ are action Lie algebroids
corresponding to actions which define representations $\rho$, of $A$
on $B$, and $\sigma$, of $B$ on $A$, with respect to which $A$ and $B$
form a matched pair.

Conversely, let $A$ and $B$ be a matched pair of Lie algebroids with
respect to $\rho$ and $\sigma$. Then the action of $A$ on $q_B$ induced
by $\rho$ and the action of $B$ on $q_A$ induced by $\sigma$ define Lie
algebroid structures on $\mathcal{A} = A\pback{M}B$ with respect to which
$(\mathcal{A};A,B;M)$ is a vacant double Lie algebroid.
\end{thm}

\pf
It remains to prove the second statement. We first verify the bialgebroid
condition (\ref{eq:full}).

For (\ref{eq:full}$)_0$ it suffices to retrace the steps of the direct
argument. Now (\ref{eq:full}) with $\phi_1 = \phi_2 = 0$ and any
$(X_1 + \psi_1, X_2 + \psi_2) = (X_1 + 0, 0 + \psi_2) +
(0 + \psi_1, X_2 + 0)$ follows by skew--symmetry of all terms in
(\ref{eq:full}).

For (\ref{eq:full}) with $Y_1 = Y_2 = 0$ and any arguments, it is
straightforward to verify that all terms are identically zero.

Consider the case $\phi_1 = 0,\ Y_2 = 0$, with arbitrary arguments.
Calculating by the same methods as in the proof of the direct statement,
and expressing all contragredient expressions in terms of $\rho$ and
$\sigma$, we see that we must prove that
\begin{equation*}
\begin{split}
-a(X_1)b(Y_1)\langle\phi_2, X_2\rangle
   & + a(X_2)b(Y_1)\langle\phi_2, X_1\rangle
          +a(X_1)\langle\phi_2, \sigma_{Y_1}(X_2)\rangle\\
   & -a(X_2)\langle\phi_2, \sigma_{Y_1}(X_1)\rangle
         +b(Y_1)\langle\phi_2, [X_1, X_2]\rangle
            -\langle\phi_2, \sigma_{Y_1}[X_1, X_2]\rangle\\
=   -b(Y_1)a(X_1) & \langle\phi_2, X_2\rangle
     +b(Y_1)a(X_2)\langle\phi_2, X_1\rangle
            +b(Y_1)\langle\phi_2, [X_1, X_2]\rangle\\
   &      \quad +a(\sigma_{Y_1}(X_1))\langle\phi_2, X_2\rangle
            -a(X_2)\langle\phi_2, \sigma_{Y_1}(X_1)\rangle
            -\langle\phi_2, [\sigma_{Y_1}(X_1), X_2]\rangle \\
   &      \quad +a(X_1)\langle\phi_2, \sigma_{Y_1}(X_2)\rangle
            -a(\sigma_{Y_1}(X_2))\langle\phi_2, X_1\rangle
            -\langle\phi_2, [X_1, \sigma_{Y_1}(X_2)]\rangle\\
   &      \quad +b(\rho_{X_2}(Y_1))\langle\phi_2, X_1\rangle
            -\langle\phi_2, \sigma_{\rho_{X_2}(Y_1)}(X_1)\rangle
            -b(\rho_{X_1}(Y_1))\langle\phi_2, X_2\rangle\\
   &      \qquad +\langle\phi_2, \sigma_{\rho_{X_1}(Y_1)}(X_2)\rangle.
\end{split}
\end{equation*}
Here six terms cancel in pairs, a further five on account of the second
equation in \ref{df:mokri}, and the remaining eight by a double application
of the last equation in \ref{df:mokri}. This completes the proof that
$(A\sdp B^*, \Bar{A^*\pds B})$ is a Lie bialgebroid.

The remaining parts of Condition~III, and Condition~I, are straightforward.
We verify Condition~II. We need the following result.

\begin{prop}{\bf \cite[2.5]{MackenzieX:1998}}
Let $(\xi, x)$ be a linear vector field on a vector bundle $E\to M$, and
let $D$ be the corresponding element of $\Ga\CDO(E)$. Then for all
$\mu\in\Ga E, m\in M$,
\begin{equation}                               \label{eq:swapE}
\xi(\mu(m)) = T(\mu)(x(m)) - D(\mu)^\upa(\mu(m)).
\end{equation}
\end{prop}

In particular for any vector fields $x, y$ on $M$, and the complete lift
$\Tilde{x} = J\circ T(x)\in\mathcal{X}(TM)$, where $J\colon T^2M\to T^2M$
is the canonical involution,
\begin{equation}                                \label{eq:swapTM}
T(y)(x(m)) = \Tilde{x}(y(m)) + [x,y]^\upa(y(m)).
\end{equation}
Because $T(y)(x(m))$ and $\Tilde{x}(y(m))$, as elements of $T^2M$, both
project to $x(m)$ under $T(p)$ and to $y(m)$ under $p_{TM}$, it follows
from (\ref{eq:swapTM}) and the interchange law for $T^2M$ that
\begin{equation}                                \label{eq:swapTMJ}
T(y)(x(m)) = \Tilde{x}(y(m)) \dpl [x,y]\hatt(x(m)),
\end{equation}
where for any vector bundle $E\to M$ and $\mu\in\Ga E$, $\Hat{\mu}$
denotes the section of $TE\to TM$ with $\Hat{\mu}(z) =
T(0)(z) + \Bar{\mu(m)}$ for $z\in T_mM$; see \cite[\S2]{MackenzieX:1998}.

Returning to the proof of Theorem~\ref{theorem:mp}, we must show that the
anchor $\tila\colon q_A\pds B \to TB$ is a Lie algebroid morphism over
$a\colon A\to TM$, where the domain is the action Lie algebroid for the
action of $B$ on $q_A$, and the target is the prolongation of $B$. First we
verify that $\tila$ commutes with the anchor $\tilb$ of the domain and the
anchor $J\circ T(b)$ of the target. Applying (\ref{eq:swapE}) to
(\ref{eq:rho}) we obtain, for any $X\in\Ga A,\ Y\in\Ga B$,
\begin{equation}                                \label{eq:tila}
\tila(X(m), Y(m)) = \xi(X)(Y(m))
   = T(Y)(aX(m)) - \rho_X(Y)^\upa(Y(m)).
\end{equation}
Then applying $T(b)$ to both sides and writing $x = aX, y = bY$, we have
$$
(T(b)\circ\tila)(X(m), Y(m))
   = T(y)(x(m)) - (b(\rho_X(Y)))^\upa(y(m))
$$
where we used the fact that $T(b)$ is a morphism of double vector bundles
$TB\to T^2M$ with core $b$.

Using the third equation of \ref{df:mokri}, and (\ref{eq:swapTM}),
this becomes
$$
\tilx(y(m)) - (a(\sigma_Y(X)))^\upa(y(m)) =
   \tilx(y(m)) \dminus (a(\sigma_Y(X)))\hatt(x(m)),
$$
where $\dminus$ is the prolonged subtraction in $T(p)\colon T^2M\to TM$.
Applying $J$ to this, then applying (\ref{eq:swapE}) to (\ref{eq:sigma})
and proceeding as before, gives
$$
T(x)(y(m)) - (a(\sigma_Y(X)))^\upa(x(m)) = T(a)(\eta(Y)(X(m)))
= T(a)(\tilb(X(m), Y(m))).
$$
Thus $T(a)\circ\tilb = J\circ T(b)\circ\tila$ as required.

We now verify the bracket condition, as given in \cite{HigginsM:1990a}.
Take $Y\in\Ga B$ and denote by $\Bar{Y}$ the pullback section of
$q_A\pds B\to A$. (It suffices to consider such sections, since they
generate $\Ga_A(q_A\pds B)$.) In general there is no section of $TB\to TM$
to which $\Bar{Y}$ projects under $\tila$. However from (\ref{eq:tila})
it follows that
$$
\tila\circ\Bar{Y} = 1\otimes T(Y)\dminus R_Y
$$
where $1\otimes T(Y)$ is the pullback of $T(Y)\in\Ga_{TM}TB$ across $a$,
and $R_Y\in\Ga(a^! TB)$ is defined by
$$
R_Y(X) = (X, \rho_X(Y)\hatt(aX)).
$$
As with any section of a pullback, $R_Y$ has a tensor decomposition
which here can be taken to be of the form
$$
\sum \ell_{\phi_i}\otimes \Hat{W_i} \dpl
\sum (f_j\circ q)\otimes \Hat{V_j}
$$
where $\phi_i\in\Ga A^*,\ f_j\in C(M),\ W_i, V_j\in\Ga B$ and
$\ell_\phi$ denotes the linear function
$A\to\R,\ X\mapsto\langle\phi(qX), X\rangle.$
(The addition $\dpl$ is the prolonged addition in $TB\to TM.$)
Now $(f\circ q)\otimes \Hat{V}
= (f\circ p \circ a)\otimes \Hat{V} = 1\otimes (f\circ p)\dtimes \Hat{V}
= 1\otimes \Hat{(fV)}$, where $p\colon TM\to M$ is the projection. The
second group of terms can therefore be collapsed to a single
$1\otimes \Hat{V}$. Since $R(0_m) = (0,0)$ and $\ell_\phi(0_m) = 0$
for all $0_m\in A$, it follows that $V = 0$ and we can take, for
$Y_1,\ Y_2\in\Ga B$,
\begin{equation}                             \label{eq:RR}
R_{Y_1} = \sum\ell_{\phi_i}\otimes \Hat{W_i},\qquad
R_{Y_2} = \sum\ell_{\psi_j}\otimes \Hat{Z_j}.
\end{equation}

Following \cite[1.3]{HigginsM:1990a}, we must prove that
$\tila\circ\Bar{[Y_1, Y_2]}$ is given by
\begin{gather}
\begin{split}\nonumber
1\otimes[T(Y_1), T(Y_2)]
& \dminus \sum\ell_{\phi_i}\otimes [\Hat{W_i}, T(Y_2)]
  \dminus \sum\ell_{\psi_j}\otimes [T(Y_1), \Hat{Z_j}]\\
& \quad \dpl \sum\ell_{\phi_i}\ell_{\psi_j}\otimes [\Hat{W_i}, \Hat{Z_j}]
  \dminus \sum \tilb(\Bar{Y_1})(\ell_{\psi_j})\otimes \Hat{Z_j}
  \dpl \sum \tilb(\Bar{Y_2})(\ell_{\phi_i})\otimes \Hat{W_i}
\end{split}\\
\begin{split}
= 1\otimes T([Y_1, Y_2])
& \dpl \sum \ell_{\phi_i}\otimes[Y_2, W_i]\hatt
\dminus \sum \ell_{\psi_j}\otimes[Y_1, Z_j]\hatt  \label{eq:12way} \nonumber \\
& \dminus \sum \eta(Y_1)(\ell_{\psi_j})\otimes \Hat{Z_j}
\dpl \sum \eta(Y_2)(\ell_{\phi_i})\otimes \Hat{W_i}
\end{split}
\end{gather}
where we used the definition of $\eta$ in terms of the anchor $\tilb$
and the following relations, from \cite[(27)]{MackenzieX:1994}, for the
bracket in $TB\to TM$:
$$
[T(Y_1), T(Y_2)] = T([Y_1, Y_2]),\quad
[T(Y_1), \Hat{Y_2}] = [Y_1, Y_2]\hatt,\quad
[\Hat{Y_1}, \Hat{Y_2}] = 0,
$$
where $Y_1,\ Y_2\in\Ga B$. Now equations (\ref{eq:RR}) are equivalent to
$$
\rho_X(Y_1) = \sum \langle \phi_i, X\rangle W_i,\qquad
\rho_X(Y_2) = \sum \langle \psi_j, X\rangle Z_j,
$$
where $X\in\Ga A,\ Y_1, Y_2\in\Ga B$, so it follows that
\begin{gather*}
[\rho_X(Y_1), Y_2] = \sum\langle \phi_i, X\rangle [W_i, Y_2]
- \sum b(Y_2)\langle\phi_i, X\rangle W_i,\\
\rho_{\sigma_{Y_2}(X)}(Y_1) =
\sum \langle\phi_i, \sigma_{Y_2}(X)\rangle W_i,
\end{gather*}
together with two similar equations. Using the first equation of
\ref{df:mokri} this leads to
\begin{multline*}
\rho_X([Y_1, Y_2]) =
         -\sum\langle\phi_i, X\rangle[Y_2, W_i]
          -\sum b(Y_2)\langle\phi_i, X\rangle W_i
          +\sum\langle\psi_j, X\rangle[Y_1, Z_j] \\
          +\sum b(Y_1)\langle\psi_j, X\rangle Z_j
   +\sum \langle\phi_i, \sigma_{Y_2}(X)\rangle W_i
    -\sum \langle\psi_j, \sigma_{Y_1}(X)\rangle Z_j,
\end{multline*}
and by (\ref{eq:sigma}) this is
$$
-\sum\langle\phi_i, X\rangle[Y_2, W_i] +\sum\langle\psi_j, X\rangle[Y_1, Z_j]
                   -\sum(\eta(Y_2)(\ell_{\phi_i})\circ X) W_i
+\sum(\eta(Y_1)(\ell_{\psi_j})\circ X) Z_j.
$$
So $\tila\circ\Bar{[Y_1, Y_2]} =
1\otimes T([Y_1, Y_2]) \dminus R_{[Y_1, Y_2]}$
indeed coincides with the RHS of (\ref{eq:12way}),
and brackets are preserved.

This completes the proof of Theorem~\ref{theorem:mp}.
\boom

In the process we have also proved the following result.

\begin{thm}                                              \label{theorem:sdp}
Let $A$ and $B$ be Lie algebroids on the same base $M$, and let
$\rho\colon A\to\CDO(B_0)$ and $\sigma\colon B\to\CDO(A_0)$ be
representations, where $A_0$ and $B_0$ denote the underlying vector
bundles.

Then $A$ and $B$ form a matched pair with respect to $\rho$ and $\sigma$
if and only if $(A\sdp B^*, \Bar{A^*\pds B})$, with the Lie algebroid
structures described in {\rm \ref{lem:sdps}}, form a Lie bialgebroid.
\end{thm}

Thus the three matched pair equations in \ref{df:mokri} are embodied in the
bialgebroid equation (\ref{eq:full}).

Theorem \ref{theorem:sdp} in the Lie algebra case has been found very recently
by Stachura \cite{Stachura}, and in the setting of Lie--Rinehart algebras
by Huebsch\-mann \cite{Huebschmann:9811069}.

In the case of a Lie bialgebra $({\mathfrak g}, {\mathfrak g}^*)$, the Lie
bialgebroid associated to the vacant double Lie algebroid is
$({\mathfrak g}\oplus{\mathfrak g}_0, {\mathfrak g}^*\oplus{\mathfrak g}^*_0)$
where the subscripts denote the abelianizations. This is of course consistent
with \ref{theorem:Manin} in the bialgebra case---which is both a
bialgebroid and a matched pair.

One other example which should be mentioned briefly is that of affinoids.
An {\em affinoid} \cite{Weinstein:1990} may be regarded as a vacant double
Lie groupoid in which both side groupoids $H$ and $V$ are the graphs of
simple foliations defined by surjective submersions $\pi_1\colon M\to Q_1$
and $\pi_2\colon M\to Q_2$. The corresponding double Lie algebroid was
calculated in \cite{Mackenzie:2000w} to be a pair of conjugate flat
partial connections adapted to the two foliations. The bialgebroid in
this case is $(T^{\pi_1}M\oplus(T^{\pi_2}M)^*,
(T^{\pi_1}M)^*\oplus T^{\pi_2}M)$ with semi--direct structures defined
by the connections.

Theorem \ref{theorem:mp} actually shows the equivalence of three formulations:
a vacant double Lie algebroid structure on $(\mathcal{A};A,B;M)$ is
equivalent to a matched pair structure on $(A, B)$, which in turn is
equivalent to a Lie bialgebroid structure on $(A\sdp B^*; A^*\pds B)$
of the type in \ref{lem:sdps}. This provides confirmation, we believe, of
the correctness of the definition of double Lie algebroid. It also provides
a diagrammatic characterization of matched pairs of Lie algebroids, directly
comparable to the diagrammatic description of matched pairs of group(oid)s
given in \cite[\S2]{Mackenzie:1992}. In the groupoid case, the twisted
multiplicativity equations are fairly unintuitive, and we believe that the
derivation of them directly from the vacant double groupoid axioms has been
a significant clarification. In the Lie algebroid case the equations in
\ref{df:mokri} are again not simple and, unlike the groupoid case, are
defined in terms of sections rather than elements. Nonetheless the
characterization given by \ref{theorem:mp} is purely diagrammatic:
recall that the characterization \ref{theorem:6.2} of a Lie bialgebroid is
formulated entirely in terms of the Poisson tensor and the canonical
isomorphism $R$. Thus we have a definition of matched pair which can be
formulated more generally in a category possessing pullbacks and suitable
additive structure.

\newcommand{\noopsort}[1]{} \newcommand{\singleletter}[1]{#1}

\end{document}